\documentclass[11pt,a4paper,oneside]{article}

\usepackage[dvips]{graphicx}
\usepackage[T1]{fontenc}
\usepackage{textcomp}
\usepackage[noBBpl]{mathpazo}
\usepackage{booktabs}
\usepackage{stylebib}

\usepackage{amsmath}
\usepackage{amsfonts}
\usepackage{amsthm}

\newcommand{\thefont}[2]{\fontsize{#1}{#2}\fontshape{n}\selectfont}
\newcommand{\1}{\rlap{\thefont{10pt}{12pt}1}\kern.16em\rlap{\thefont{11pt}{13.2pt}1}\kern.4em}

\renewcommand{\baselinestretch} {1.5}

\makeatletter
\def\singlespace{\def\baselinestretch{1}\@normalsize}

\newcommand{\argmin}{\displaystyle \mathop{argmin}}
\newcommand{\tend}{\displaystyle \mathop{\longrightarrow}}
\newcommand{\somme}{\displaystyle \sum}

\renewcommand{\epsilon}{\varepsilon} 
\renewcommand{\phi}{\varphi}

\def\1{\mbox{1\hspace{-.35em}1}}
\def\R{\mathbb{R}}

\def\N{\mathbb{N}}
\def\P{\mathbb{P}}
\def\E{\mathbb{E}}
\def\L{\mathbb{L}}

\def\R{\mathbb{R}}

\newcommand{\ints}{\ensuremath{{\mathbb Z}}}

\newcommand{\bd}{\mbox{${\mathbf d}$}}

\newcommand{\be}{\mbox{${\mathbf e}$}}

\newcommand{\bz}{\mbox{${\mathbf z}$}}

\newcommand{\bA}{\mbox{${\mathbf A}$}}

\newcommand{\bR}{\mbox{${\mathbf R}$}}

\newcommand{\bU}{\mbox{${\mathbf U}$}}

\newcommand{\bff}{\mbox{${\mathbf F}$}}

\newcommand{\bY}{\mbox{${\mathbf Y}$}}

\newcommand{\bX}{\mbox{${\mathbf X}$}}

\newcommand{\bZ}{\mbox{${\mathbf Z}$}}

\newcommand{\balpha}{\mbox{\boldmath$\alpha$}}
\newcommand{\bbeta}{\mbox{\boldmath$\beta$}}

\newcommand{\btheta}{\mbox{\boldmath$\theta$}}

\newcommand{\bepsilon}{\mbox{\boldmath$\epsilon$}}

\newcommand{\la}{\mbox{$\langle$}}

\newcommand{\ra}{\mbox{$\rangle$}}

\newtheorem{Theoreme}{Theorem}
\newtheorem{Proposition}{Proposition}[]

\newtheorem{Lemme}[Proposition]{Lemma}

\newtheorem{rmq}{Remark}

\title{{\sc Robust estimation and Wavelet Thresholding in Partial Linear Models}}

\author{Ir\`ene Gannaz\\ 
{\em Laboratoire de Mod\'elisation et Calcul}\\
{\em Universit\'e Joseph Fourier}\\
{\em BP 53 - 38041 Grenoble Cedex 9}\\ {\em France}} 

\date{November 2006}

\newcommand{\figs}[3]{\begin{figure}[htbp] %
 \begin{center}\includegraphics[width=#2]{#1.eps}\end{center} %
 \small %
 \begin{singlespace} %
 \caption{#3}
 \label{#1} %
 \end{singlespace} %
 \end{figure}}

\newcommand{\twofig}[4]
{ \hbox to\hsize{\hss
 \vbox{\psfig{figure=#1,width=#3,height=#4}}\qquad
 \vbox{\psfig{figure=#2,width=#3,height=#4}}
 \hss}
\vskip -0.0truein \hbox to\hsize{\hss
 \vbox{ \begin{center}\mbox{\footnotesize \hspace{0.0in} {(a)}
 \hspace{#3} {(b)} } \end{center} }
 \hss}
\vskip 0.0truein }
\newcommand{\threefig}[5]{
\hbox to\hsize{\hss
 \vbox{\psfig{figure=#1,width=#4,height=#5}} \hspace{0.0in}
 \vbox{\psfig{figure=#2,width=#4,height=#5}} \hspace{0.2in}
 \vbox{\psfig{figure=#3,width=#4,height=#5}}
 \hss}
\vskip -0.0truein \hbox to\hsize{\hss
 \vbox{ \begin{center}\mbox{\footnotesize \hspace{0.0in} {(a)}
 \hspace{#4} {(b)} \hspace{#4} {(c)}} \end{center} }
 \hss}
}
\newcommand{\fourfig}[6]
{ \hbox to\hsize{\hss
 \vbox{\psfig{figure=#1,width=#5,height=#6}}\qquad
 \vbox{\psfig{figure=#2,width=#5,height=#6}}
 \hss}
\vskip -0.0truein \hbox to\hsize{\hss
 \vbox{ \begin{center}\mbox{\footnotesize \hspace{0.1in} {(a)}
 \hspace{#5} {(b)} } \end{center} }
 \hss}
\vskip 0.1truein \hbox to\hsize{\hss
 \vbox{\psfig{figure=#3,width=#5,height=#6}}\qquad
 \vbox{\psfig{figure=#4,width=#5,height=#6}}
 \hss}
\vskip -0.1truein
 \vbox{ \begin{center}\mbox{\footnotesize \hspace{0.1in} {(c)}
 \hspace{#5} {(d)} } \end{center} }
\hbox to\hsize{\hss
 \hss}
\vskip -0.1truein }

\begin{document}

\maketitle

\begin{abstract}
{\small This paper is concerned with a semiparametric partially linear regression 
model with unknown regression coefficients, an unknown nonparametric function 
for the non-linear component, and unobservable Gaussian distributed random errors. 
We present a wavelet thresholding based estimation 
procedure to estimate the components of the partial linear 
model by establishing a connection between an $l_1$-penalty based wavelet 
estimator of the nonparametric component and Huber's M-estimation of a standard linear model with outliers. Some general results on the large sample properties of the estimates of both the parametric and the nonparametric part of the model are established. 
Simulations and a real example are used to illustrate the general results and to compare the proposed
methodology with other methods available in the recent literature.}
\end{abstract}

\noindent {\it Keywords}: Semi-nonparametric models, partly linear models, wavelet thresholding, backfitting, M-estimation, penalized least-squares. 

\medskip

\section{Introduction} 

Assume that responses $y_1,\ldots,y_n$ are observed at deterministic equidistant points $t_i=\frac i n$ of an univariate variable such as time and for fixed values $\mathbf{X}_i$, $i=1, \ldots, n$, of some $p$-dimensional explanatory variable and that the relation between the response and predictor values is modeled by a Partially Linear Model (PLM): 
\begin{equation}\label{plm}
y_i=\mathbf{X}_i^T\bbeta_0+f(t_i)+u_i\qquad i=1,\ldots, n,
\end{equation}
where $\bbeta_0$ is an unknown $p$-dimensional real
parameter vector and $f(\cdot)$ is an unknown real-valued function; the $u_i$'s are i.i.d. normal errors with mean 0 and variance $\sigma^2$ and superscript ``T'' denotes the transpose of a vector or matrix. 
Given the observed data $(y_i,\mathbf{X}_i)_{i=1,\ldots, n}$, the aim is to estimate from the data the vector $\bbeta$ and the function $f$.

The interest in partial linear models has grown significantly within the last decade since their introduction by \shortciteA{Engle} to analyze in a nonlinear fashion the relation between electricity usage and average daily temperature. Since then the models have been widely studied in the 
literature. The recent monograph by \shortciteA{Hardle} provides an excellent survey on the theory and applications of the model in a large variety of fields, such as finance, economics, geology and biology, to name only a few. The advantages of such a model is that it allows
an adequate and more flexible handling of the explanatory variables than in linear models and can be 
also serve as a starting point for dimension reduction by additive modeling. Although there is still lack of general theory on testing the goodness-of-fit of a partial linear model, there are some 
consistent specification tests such as, for example, those developed by \shortciteA{ChenChen}.

Until now, several methods have been proposed to analyse partially linear models.
One approach to estimation of the nonparametric component in these models is based on smoothing splines regression techniques and has been employed in particular by \shortciteA{GreenYandell}, \shortciteA{Engle}, \shortciteA{Rice}, \shortciteA{Chen87}, \shortciteA{ChenShiau}, and \shortciteA{Schick} among others. Kernel regression (see e.g. \shortciteA{Speckman}) and local polynomial fitting techniques (see e.g. \shortciteA{HamiltonTruong}) have also been used to study partially linear models. An important assumption by all these methods for the unknown nonparametric component $f(t)$ is its high smoothness. But in reality, such a strong assumption may not be satisfied. 
To deal with cases of a less-smooth nonparametric component, a wavelet based estimation procedure
is developed in this paper, and as such it can handle nonparametric estimation for curves lying in Besov spaces instead of the more classical Sobolev spaces.

The estimation method developed in this paper is based on a wavelet expansion of the nonparametric part of
the model. The use of an appropriate thresholding strategy on the coefficients allows us to estimate in an adaptive way the nonparametric part with quasi-minimax asymptotic rates without restrictive assumptions on its regularity. To our knowledge, only few developments in the use of nonlinear wavelet 
methods in the context of PLM models exist in the literature. Wavelet based estimators for the nonparametric 
component of a PLM have been investigated by \shortciteA{Meyer03}, \shortciteA{ChangQu} and by \shortciteA{FadiliBullmore}, more recently. Our results will be compared to the later, since the settings adopted in their work are relatively similar to ours.

One novelty of the estimation procedure proposed in this paper is the link between wavelet thresholding and classical robust M-estimation schemes in linear models with outliers: using soft or hard thresholding or even a SCAD thresholding (see \shortciteA{AntoniadisFan}) amounts in estimating respectively the unknown vector $\bbeta_0$ of the linear part in the model by Huber's M-estimation or by a truncated mean or by Hampel's estimator. This link allows us to investigate the asymptotic minimax properties of the estimators and to derive second-order approximations for the bias and variance of the resulting estimators of $\bbeta_0$. This is essentially due to the fact that the nonparametric part of the model has a sparse wavelet coefficients representation, and the wavelet coefficients of a PLM in the wavelet domain appear then as outliers in the linear model composed by the linear part. 

Furthermore, the above established link of our method with M-estimation theory offers the possibility to use specific M-estimation algorithms for numerically implementing the proposed method, instead of using the {\em backfitting} technique proposed by \shortciteA{FadiliBullmore}. For our numerical implementation, we will adopt a class of half-quadratic optimization algorithms that have been developed recently for robust image recognition in the pattern recognition literature (see e.g. \shortciteA{Charbonnier}, \shortciteA{Dahyot-comp}, \shortciteA{Vik} and \shortciteA{NikolovaNg}).

The organization of this paper is as follows: Section 2 briefly recalls
some relevant facts about the wavelet series expansion and the
discrete wavelet transform that we need further and presents the wavelet decomposition
used to model the observed partial linear model. In section 3, we establish the connection between wavelet thresholding estimation for the PLM and M-estimation for a linear model. Section 4 establishes the main properties of our estimators. In Section 5, we discuss the computational algorithms that are used for the numerical implementation of our procedures where we also present a small simulation study to illustrate the finite sample properties of our procedures and to compare them to the backfitting algorithm proposed by \shortciteA{FadiliBullmore}. Proofs of our results are given in Appendix. 

\section{The partly linear model and its wavelet transform}

\subsection{{\sc The Setup}}
\label{subsec:s}

Suppose that $y_{i}$ ($i=1,2,\ldots,n$) is
the $i$-th response of the regression model at point $t_{i}$ (where $t$ is
an index such as time or distance) and can be modelled as
\begin{equation} \label{PLM}
y_{i} = \bX_{i}^T \bbeta_0+ f(t_{i}) + u_i,
\end{equation}
where $ \bX_{i}^T$ are given $p \times 1$ vectors of covariate values, $t_i =\frac i n$
 and $\bbeta_0$ and $f$ are respectively
the parametric and nonparametric components of the partial linear model. We will assume hereafter that
 the noise variables $u_i$ are i.i.d. Gaussian $\mathcal N(0,\sigma^2)$ and that the sample size
$n=2^J$ for some positive integer $J$.

In the nonparametric analysis, the nonparametric part $f$ is modeled as a function lying in
an infinite dimensional space. The underlying notion behind wavelet methods is that the unknown function 
has an economical wavelet expression, i.e. $f$ is, or is well approximated by, 
a function with a relatively small proportion of nonzero wavelet coefficients. 
An approach to modelling the nonparametric component of the PLM model, that allows a wide range
of irregular effects, is through the sequence space representation of Besov spaces. The (inhomogeneous) Besov spaces on the unit interval, $\mathcal B^s_{\pi,r}([0,1])$, consist of functions
that have a specific degree of smoothness in their derivatives. The parameter
$\pi$ can be viewed as a degree of function's inhomogeneity while $s$ is a
measure of its smoothness. Roughly speaking, the (not necessarily integer)
parameter $s$ indicates the number of function's (fractional) derivatives,
where their existence is required in an $L^{\pi}$-sense; the additional
parameter $r$ is secondary in its role, allowing for additional fine
tuning of the definition of the space. For a detailed study on (inhomogeneous)
Besov spaces we refer to, e.g., \shortciteA{DonJohn98}. To capture key
characteristics of variations in $f$ and to exploit its sparse wavelet coefficients
representation, we will assume that $f$ belongs to $\mathcal B^s_{\pi,r}([0,1])$ with $s+1/\pi-1/2>0$.
The last condition ensures in particular that evaluation of $f$ at a given point makes sense.

We now briefly recall first some relevant facts about the wavelet series expansion and the discrete wavelet transform that we need further.

\subsection{The wavelet series expansion and the discrete wavelet transform}
\label{subsec:wavelets}

Throughout the paper we assume that we are working within an orthonormal basis
generated by dilatations and translations of a compactly supported scaling
function, $\phi(t)$, and a compactly supported mother wavelet, $\psi(t)$,
associated with an $r$-regular ($r \geq 0$) multiresolution analysis of
$\left(L^{2}[0,1], \la \cdot, \cdot \ra\right)$, the space of
squared-integrable functions on $[0,1]$ endowed with the inner product $\la f,
g \ra = \int_{[0,1]}f(t)g(t)\,dt$. For simplicity in exposition, we work with
periodic wavelet bases on $[0,1]$ (see, e.g., \shortciteA{Mallat99}, Section 7.5.1),
letting
$$
\phi_{jk}^{{\rm p}}(t) = \sum_{l \in \ints} \phi_{jk}(t-l) \quad
\text{and} \quad \psi_{jk}^{{\rm p}}(t) = \sum_{l \in \ints} \psi_{jk}(t-l),
\quad \mbox{for} \quad t \in[0,1],
$$
where $ \phi_{jk}(t) = 2^{j/2}\phi(2^{j}t-k)$ and $\psi_{jk}(t) =
2^{j/2}\psi(2^{j}t-k)$. For any given primary resolution level $j_{0} \geq 0$,
the collection
$$
\{\phi_{j_{0}k}^{{\rm p}},\; k=0,1,\ldots,2^{j_{0}}-1;\; \psi_{jk}^{{\rm p}},\; j \geq j_{0};
\; k=0,1,\ldots,2^{j}-1 \}
$$
is then an orthonormal basis of $L^{2}[0,1]$. The superscript
``${\rm p}$'' will be suppressed from the notation for
convenience. Despite the poor behavior of periodic wavelets near
the boundaries, where they create high amplitude wavelet
coefficients, they are commonly used because the numerical
implementation is particular simple. Therefore, for any $f \in
L^2[0,1]$, we denote by $c_{j_0k} = \la f, \phi_{j_0k} \ra $
($k=0,1,\ldots,2^{j_0}-1$) the scaling coefficients and by $d_{jk}
= \la f, \psi_{jk} \ra $ ($j \geq j_0$;\;
$k=0,1,\ldots,2^{j}-1$) the wavelet coefficients of $f$ for the
orthonormal periodic wavelet basis defined above; the function
$f$ is then expressed in the form
$$
f(t) = \sum_{k=0}^{2^{j_0}-1} c_{j_0k}\phi_{j_0k}(t) +
\sum_{j=j_0}^{\infty} \sum_{k=0}^{2^{j}-1} d_{jk}\psi_{jk}(t), \quad t \in[0,1].
$$
The approximation space spanned by the scaling functions $\{\phi_{j_{0}k}, \; k=0,1,\ldots,2^{j_{0}}-1\}$ is usually denoted by $V_{j_0}$ while the details space at scale $j$, spanned by $\{ \psi_{jk},\;
\; k=0,1,\ldots,2^{j}-1 \}$ is usually denote by $W_j$. 

In statistical settings, we are more usually concerned with discretely sampled,
rather than continuous, functions. It is then the wavelet analogy to the
discrete Fourier transform which is of primary interest and this is referred to
as the discrete wavelet transform (DWT). Given a vector of real values
$\be=(e_1,\ldots, e_n)^T$, the discrete
wavelet transform of $\be$ is given by $\bd=W_{n \times n} \be$, where $\bd$
is an $n \times 1$ vector comprising both discrete scaling coefficients,
$s_{j_0k}$, and discrete wavelet coefficients, $w_{jk}$, and $W_{n \times n}$
is an orthogonal $n \times n$ matrix associated with the orthonormal periodic
wavelet basis chosen. In the following we will distinguish the blocs of $W_{n \times n}$ spanned
respectively by the scaling functions and the wavelets.
The empirical coefficients $s_{j_0k}$ and $w_{jk}$ of $\be$ are 
given by
\begin{eqnarray*}
s_{j_0,k}&\approx&\frac{1}{\sqrt{n}} \somme_{i=1}^n e_i\phi_{j_0,k}(t_i) \quad \text{for}~~k=0, \ldots, 2^{j_0}-1\\
w_{j,k}&\approx &\frac{1}{\sqrt{n}} \somme_{i=1}^n e_i\psi_{j,k}(t_i) \quad \text{for}~~\left\{\begin{array}{rcl}
j&=&j_0, \ldots, J-1,\\
k&=&0, \ldots, 2^{j}-1.
\end{array}\right.\\
\end{eqnarray*} 
When $\be$ is a vector of function values $\bff=(f(t_1),...,f(t_n))^T$ at equally spaced points $t_i$, 
the corresponding empirical coefficients $s_{j_0k}$ and $w_{jk}$ are
related to their continuous counterparts $c_{j_0k}$ and $d_{jk}$ (with an approximation error of
order $n^{-1}$) via the relationships $s_{j_0k} \approx \sqrt{n}\, c_{j_0k}$
and $w_{jk} \approx \sqrt{n}\, d_{jk}$. Note that, because of orthogonality of
$W_{n \times n}$, the inverse DWT (IDWT) is simply given by $\bff=W_{n \times
n}^{\rm T} \bd$, where $W_{n \times n}^{\rm T}$ denotes the transpose of $W_{n
\times n}$. If $n=2^J$ for some positive integer $J$, the DWT and IDWT may be
performed through a computationally fast algorithm (see, e.g., \shortciteA{Mallat99},
Section 7.3.1) that requires only order $n$ operations.

We will further use the following notation. For a $n$-dimensional vector $\be$, its Euclidian (or $l_2$) norm $\left(\sum_{i=1}^n e_i^2 \right)^{1/2}$ will be denoted by $\| \be \|$ and the Frobenius norm of a matrix $B$ with general entries $b_{i,j}$ will be denoted by $\|B\|=\left(\sum_{i,j} b_{i,j}^2 \right)^{1/2}$.

\subsection{A wavelet-based model specification of the PLM model}
\label{subsec:wbms}
In matrix notation, the PLM model specified by (\ref{PLM}) can be written as
\begin{equation}
\label{modele0}
\bY=X\bbeta_0+\bff+\bU,
\end{equation}
where $\bY=\begin{pmatrix} y_1, \ldots, y_n\end{pmatrix}^T$, $X^T=\begin{pmatrix} \bX_1, \ldots, \bX_n\end{pmatrix}$ is the $p\times n$ design matrix, and $\bff=\begin{pmatrix} f(t_1), \ldots f(t_n)\end{pmatrix}^T$. The noise vector $\bU=\begin{pmatrix} u_1, \ldots, u_n\end{pmatrix}^T$ is a Gaussian vector with mean 0 and variance matrix $\sigma^2 I_n$. 


For the model to be asymptotically identifiable, we will assume:
\begin{description}
\item[(A1)] The vector $\frac 1 n X^T\bff$ tends to 0 as $n$ goes to infinity.
\item[(A2)] The matrix $X$ is full rank, i.e. $\frac 1 n X^T X$ converges towards an invertible matrix.
\end{description}
Expressing the vector of coefficients of the linear part as
$$\bbeta_0=\left(\frac 1 n X^TX\right)^{-1}X^T(\bY-\bff-\bU),$$ 
clearly shows that conditions (A1) and (A2) are sufficient to asymptotically ensure the identifiability of the
PLM model. As it will be seen in the Appendix, none of these assumptions is restrictive.

Let now $\bZ=W_{n \times n}\bY$, $A=W_{n \times n}X$, $\btheta_0=W_{n \times n}\bff$ and $\bepsilon= W_{n \times n}\bU$. Then premultiplying (\ref{plm}) by $W$, we obtain the transformed model 
\begin{equation}
\label{Wmodele} 
\bZ=A\bbeta_0+\btheta_0+\bepsilon.
\end{equation}
The orthogonality of the DWT matrix $W_{n \times n}$ ensures that the transformed noise vector $\bepsilon$
is still distributed as a Gaussian white noise with variance $\sigma^2 I_n$. Hence, the representation of the model in the wavelet domain not only allows to retain the partly linear structure of the model but also to exploit
in an efficient way the sparsity of the wavelet coefficients in the representation of the nonparametric component. 

\section{Soft Thresholding and Huber's M-estimation}

The wavelet shrinkage estimators that are classically obtained by hard or soft thresholding can be regarded as an extension of the penalized least squares (PLS) estimator (see \shortciteA{AntoniadisFan}). We therefore propose estimating the parameters $\bbeta_0$ and $\btheta_0$ in model (\ref{Wmodele}) by penalized least squares. To be specific, our wavelet based estimators will be defined as follows:
\begin{equation}
\label{critere}
(\hat\bbeta_n,\hat\btheta_n)=\argmin_{(\bbeta,\btheta)}\left\{\;J_n(\bbeta,\btheta)=\sum_{i=1}^n \frac 1 2 (z_i-\bA_i^T\bbeta-\theta_i)^2+ \lambda\sum_{i=i_0}^n |\theta_{i}|\;\right\},
\end{equation} 
for a given penalty parameter $\lambda$, where $i_0=2^{j_0}+1$. The penalty term in the above expression penalizes only the empirical wavelet coefficients of the nonparametric part of the model and not its scaling coefficients. The choice $l^1$ of the penalty function produces the soft thresholding rule. 

The regularization method proposed above is closely related to the method proposed recently by 
\shortciteA{ChangQu}, but these authors essentially concentrate on the backfitting algorithms involved in the optimization, without any theoretical study of the resulting estimates. The method also relates to the recent one developed by \shortciteA{FadiliBullmore} where a variety of penalties is discussed. Note, however, that their study is limited to quadratic penalties which amounts essentially in assuming that the underlying function $f$ belongs to some Sobolev space and does not exploit the sparse representation of $f$. 

In order to establish the link with Huber's estimation we will have a closer look at the minimization of the criterion $J_n$ stated in~(\ref{critere}). For a fixed value of $\bbeta$, the criterion $J_n(\bbeta, \cdot)$ is minimum
at
\begin{equation}\label{tildetheta}\tilde \theta_i(\bbeta)=\begin{cases} z_i-\bA_i^T\bbeta & \text{if~} i<i_0,\\ \text{sign}(z_i-\bA_i^T\bbeta)\left(|z_i-\bA_i^T\bbeta|-\lambda\right)_+ & \text{if~} i\geq i_0. \end{cases}
\end{equation}
Therefore, finding $\hat\bbeta_n$, a solution to problem (\ref{critere}), amounts in finding $\hat\bbeta_n$ minimizing the criterion $J_n(\tilde\btheta(\bbeta),\bbeta)$. However, note that
\begin{equation}\label{point1} 
J_n(\tilde\btheta(\bbeta),\bbeta) =\sum_{i=i_0}^{n} \rho_\lambda (z_i-\bA_i^T\bbeta)
\end{equation} 
where $\rho_\lambda$ is Huber's cost functional with threshold $\lambda$, defined by:
\begin{equation}
\label{rho} 
\rho_\lambda(u) = \begin{cases} u^2/2 & \text{if}~~|u|\leq \lambda, \\
\lambda |u| -\lambda^2/2 & \text{if}~~|u| > \lambda. \\ \end{cases}
\end{equation}

The above facts can be derived as follows. 
Let $i\geq i_0$. Minimizing expression (\ref{critere}) with respect to $\theta_i$ is equivalent in minimizing $j(\theta_i):=\frac{1}{2}(z_i-\bA_i^T\bbeta-\theta_i)^2+\lambda |\theta_i|$. The first order condition for this is: $j'(\theta_i)=\theta_i-(z_i-\bA_i^T\bbeta)+\text{sign}(\theta_i)\lambda=0$ where $j'$ denotes the derivative of $j$. Now, 
\begin{itemize}
\item[$\bullet$] if $\theta_i\geq0$, then $j'(\theta_i)=0$ if and only if $\theta_i=z_i-\bA_i^T\bbeta-\lambda$. Hence, if $z_i-\bA_i^T\bbeta\leq\lambda$, $\theta_i=0$ and otherwise $\theta_i=z_i-\bA_i^T\bbeta-\lambda$.
\item[$\bullet$] if $\theta_i\leq0$, $j'(\theta_i)$ is zero if and only if $\theta_i=z_i-\bA_i^T\bbeta+\lambda$; therefore, if $z_i-\bA_i^T\bbeta\geq-\lambda$, $\theta_i=0$ and otherwise $\theta_i=z_i-\bA_i^T\bbeta+\lambda$.
\end{itemize}
This proves that for a fixed value of $\bbeta$, the criterion (\ref{critere}) is minimal for $\tilde\btheta(\bbeta)$ given by expression (\ref{tildetheta}).
If we now replace $\btheta$ in the 
objective function $J_n$ we obtain $J_n(\bbeta,\tilde\btheta(\bbeta))=\frac{1}{2}{\somme_{i=i_0}^{n}} \left((z_i-\bA_i^T\bbeta-\tilde\theta_i)^2+\lambda|\tilde\theta_i|\right)$ since $\tilde\theta_i=z_i-\bA_i^T\bbeta$ for $i<i_0$. Now denoting by $I$ the set $I:=\left\{j=i_0, \ldots, n, \quad |z_j-\bA_j\bbeta |<\lambda\right\}$, we find that $J_n(\bbeta,\tilde\btheta(\bbeta))=\frac{1}{2}{\somme_{I}} (z_i-\bA_i^T\bbeta)^2 + \frac{1}{2}{\somme_{I^C}} \lambda^2 +\lambda{\somme_{I^C}}\left(|z_i-\bA_i^T\bbeta|-\lambda\right)$ by replacing $\tilde\theta_i$ with (\ref{tildetheta}), which is exactly Huber's functional.

The mathematical equivalence of the solution of the two classes of estimation can be stated in the following proposition.

\begin{Proposition}\label{prop1}
If $\hat\bbeta_n$ and $\hat\btheta_n$ are solutions of the optimization problem (\ref{critere}), then they
satisfy 
\begin{eqnarray} \label{beta}
\hat\bbeta_n&=&\argmin_{\bbeta} \sum_{i=i_0}^{n} \rho_\lambda (z_i-\bA_i^T\bbeta),\\ \label{theta}
\hat\theta_{i,n}&=& \begin{cases} z_i-\bA_i^T\hat\bbeta_n & \text{if~} i<i_0\,\\ \gamma_{soft,\lambda}(z_i-\bA_i^T\hat\bbeta_n) & \text{if~} i\geq i_0, \end{cases}, \quad \; i=1,\ldots, n,
\end{eqnarray}
with $\rho_\lambda$ being Huber's cost functional defined in (\ref{rho}) and $\gamma_{soft,\lambda}$ the soft-thresholding function with threshold $\lambda$ defined by $\gamma_{soft,\lambda}(u)=\text{sign}(u)\left(|u|-\lambda\right)_+$.
\end{Proposition}
This result allows the computation of the estimators $\hat\bbeta_n$ et $\hat\btheta_n$ in a non-iterative
fashion. We can estimate the parameter $\bbeta_0$ directly from the observed data without caring about the nonparametric part of the model by means of eq.(\ref{beta}), and then determine $\hat \btheta_n$, thence $\hat \bff_n$ using eq.(\ref{theta}). 

The resulting form of the estimators allows us to study their asymptotic properties. Moreover, as we shall see
in the simulation section of this paper, another benefit is that we can design estimation algorithms that are much faster
than those based on backfitting. Lastly, Propostion~\ref{prop1} leads to a nice interpretation of the estimators.

We may summarize the estimation procedure as follows.
Using the observed data $(\bY, X)$~:
\begin{enumerate}
\item Apply the DWT of order $J=\log_2(n)$ on $X$ and $\bY$ to get their corresponding representation $A$ and $\bZ$ in the wavelet domain.
\item The parameter $\bbeta_0$ is then Huber 's robust estimator which is obtained without taking care of the nonparametric component in the PLM model, given by the optimization problem (\ref{beta}). In other words this amounts in considering the linear model $z_i=\bA_i^T\bbeta_0+e_i$ with noise $e_i=\theta_{0i}+\epsilon_i$.
\item The vector $\btheta$ of wavelet coefficients of the function $f$ is estimated by soft thresholding of $\bZ-A\hat\bbeta_n$, i.e. by equation (\ref{theta}).
The estimation of $f$ is then obtained by applying the inverse discrete wavelet transform. Note that this last step corresponds to a standard soft-thresholding nonparametric estimation of $f$ in the model: 
$$y_i-\bX_i^T\hat\bbeta_n=f(t_i)+v_i,\quad i=1, \ldots, n,$$ where $v_i=\bX_i^T(\bbeta_0-\hat\bbeta_n)+u_i$.
\end{enumerate}

\begin{rmq}
The above estimation procedure is in phase with the one advocated by \shortciteA{Speckman} who suggests that it is usually preferable to estimate first the linear component in a PLM and to then proceed to the estimation of the nonparametric one. Indeed, we propose to estimate $\bbeta_0$ and $\bff$ by: $\hat \bbeta=(X^TS^TSX)^{-1}S\bY$ and $\hat \bff= (I-S)(\bY-X\hat\bbeta)$, with $S=(I-T)W$, $T$ being the threshold operator. We recognize the exact same form as those of \shortciteA{Speckman}, differing only on the fact that the smoothing operator $S$ is not anymore linear.
\end{rmq}

The wavelet soft-thresholding procedure proposed in this section was derived by establishing the connection between an $l_1$ based penalization of the wavelet coefficients of $f$ and Huber's M-estimators in a linear model. Other penalties, leading to different thresholding procedures can also be seeing as $M$-estimation procedures. For example, if $\gamma_\lambda$ denotes the resulting thresholding function, we can show in a similar way that the estimators verify
 \begin{eqnarray*}
\hat\bbeta_n&=&\argmin_{\bbeta }\sum_{i=i_0}^{n} \rho_\lambda (z_i-\bA_i^T\bbeta),\\
\hat\theta_{i,n}&=& \begin{cases} z_i-\bA_i^T\bbeta & \text{if~} i<i_0,\\ \gamma_{\lambda}(z_i-\bA_i^T\bbeta) & \text{if~} i\geq i_0, \end{cases}, \quad i=1, \ldots, n,\\
\end{eqnarray*}
with $\rho_\lambda$ being the primitive of $u\mapsto u-\gamma_\lambda(u)$. From what precedes, one sees that
hard thresholding corresponds to mean truncation, while SCAD thresholding is associated to Hampel's
M-estimation. The above thresholding procedures and the corresponding criteria are illustrated in Figure 1.

\figs{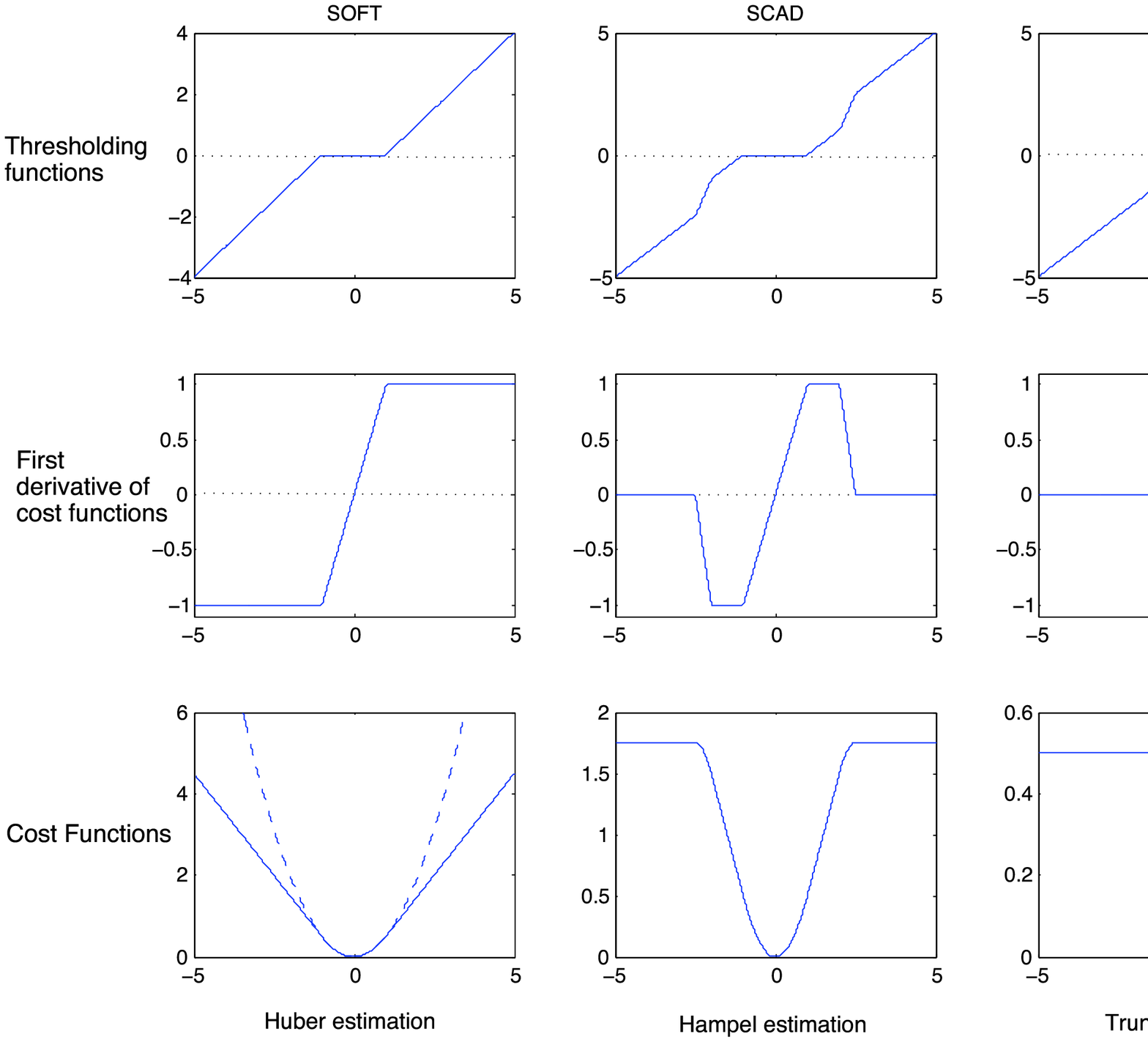}{12cm}{Link between different thresholdings and M-estimation. The dashed line displays the least squares criterion.}

However, in this paper, we only concentrate on the properties of estimators obtained by soft thresholding, 
those corresponding to other rules presenting avenues for further research that hope will be addressed in the future.


\section{Asymptotic properties}

Huber's M-estimation was introduced as an alternative to least squares in order to limit the sensitivity of the least-squares estimates to each individual observation. While Huber's M-estimators do not have finite breakdown points, one can show they are quite robust to outliers (see e.g. \shortciteA{Hampel}). 
Huber's M-estimation appears therefore a natural approach for robustly fitting the linear part of a PLM, interpreting the wavelet coefficients of the nonparametric part as outliers. In what follows, relying upon 
this analogy, we study the asymptotic properties of our estimator. To establish our asymptotic results we will require several assumptions to hold.\\

\medskip

 First a condition which ensures the unicity of $\hat\bbeta_n$ defined in (\ref{beta}):
\begin{description}
\item[(A3)] The series $(K_n)$ defined by $K_n:=\frac 1 n \sum_{i=1}^n \bA_i \bA_i^T\rho''_\lambda(\theta_{0 i}+\epsilon_i)$, converges in the $L^2$-norm towards a non-singular matrix $K_0$.
\end{description}
The next assumption deals with the structure of the regression design matrix. Since the dicrete wavelet transform $W$ is orthogonal it follows that $A^TA=X^TX$ and, therefore when (A2) holds the matrix $A^TA$ is non-singular for $n$ sufficiently large. Consequently, the projection matrix on the space spanned by the columns of $A$, say $H=A(A^TA)^{-1}A^T$, has a rank $p$. In such a case, if $(h_1, \ldots, h_n)$ denotes the diagonal of $H$, the equality $\sum h_i= p$ holds. With regards to the design,
we will also use the assumption:
\begin{description}
\item[(A4)] The quantity $h:=\displaystyle{\max_{i=1, \ldots, n}} \bA_i^T(A^TA)^{-1} \bA_i$ tends to 0 when $n$ goes to infinity.
\end{description}
Assumption (A4) is common in a robust regression framework, validating among other things the use of the Lindeberg-Feller criterion. The only difference in our case is that the regression matrix that we consider is the wavelet transformed $A$ rather than $X$, but the relevant discussion in the Appendix shows that such an assumption is reasonable.\\

\medskip

Existing results for semi-parametric partial linear models establish parametric rates of convergence for the linear part and minimax rates for the nonparametric part, showing in particular that the existence of a linear component does not changes the rates of convergence of the nonparametric component. Within the framework adopted in this paper, the rates of convergence are similar, but an extra logarithmic term will appear in the rates of the parametric part, mainly due to the fact that our smoothness assumptions on the nonparametrric part are weaker. We are now in position to give our asymptotic results. 

\begin{Theoreme}\label{resultat}
Let $\hat\bbeta_n$ and $\hat\btheta_n$ be the estimators defined by (\ref{beta},\ref{theta}) in the model (\ref{plm}). Consider that the penalty parameter $\lambda$ is the universal threshold: $\lambda=\sigma\sqrt{2\log(n)}$.
Under assumptions (A1)-- (A4), we have
\begin{eqnarray*}
\hat\bbeta_n-\bbeta_0&=&\bigcirc_\P \left(\sqrt{\frac{\log(n)}{n}}\right),\\
\text{and}\qquad\sqrt{n}(\hat\bbeta_n-\bbeta_0)&=&K_0^{-1}\left(\frac{1}{\sqrt{n}}\sum_{i=1}^n \rho_{\lambda}'(\theta_{0 i}+\epsilon_i)\bA_i\right)+o_\P(\sqrt{\log(n)}).\\
\end{eqnarray*}
If in addition we assume that the scaling function $\phi$ and the mother wavelet $\psi$ belong to $\mathcal C^R$ and that $\psi$ has $N$ vanishing moments, then, for $f$ belonging to the Besov space $\mathcal B^s_{\pi,r}$ with $0<s-1/2+1/\pi$ and $1/\pi<s<\min(R,N)$, we have
$$\|\hat f_n-f\|_2 =\bigcirc_\P\left(\left(\frac{\log(n)}{n} \right)^{\frac{s}{1+2s}}\right),$$
where $\|\hat f_n-f\|_2^2=\int_0^1 (\hat f_n-f)^2$.
\end{Theoreme}
The Theorem is proved in the Appendix.
As noted previously, we lose a factor $\sqrt{\log(n)}$ in the estimation of the vector of parameters $\bbeta$. The presence of a logarithmic loss lies on the choice of the threshold $\lambda$: taking $\lambda$ which tends to 0, as suggested by \shortciteA{FadiliBullmore}, would lead to a minimax rate in the estimation of $\bbeta$. The drawback is that the quality of the estimation for the nonparametric part of the PLM would not be anymore quasi-minimax. This phenomenon was put in evidence by \shortciteA{Rice}: a compromise must be done between the optimality of the linear part estimation with an oversmoothing of the functional estimation and a loss in the linear regression parameter convergence rate but a correct smoothing of the functional part.

The method of estimation that we propose leads to quasi-minimax convergence rates and is applicable for a large class of functions $f$. An important remark is also that our procedure is adaptative relatively to the regularity of $f$, thanks to the use of threshold techniques in the wavelet decomposition. Note also that Theorem \ref{resultat} give a Bahadur's representation of $\hat\bbeta_n$, allowing to elaborate appropriate testing procedures; such inferential problems are out of the scope of the present paper, but interesting for future work.

\subsection{Estimation of the variance}

Our estimation procedure relies upon knowledge of the variance $\sigma^2$ of the noise, appearing in the expression of the threshold $\lambda$ (recall that we have adopted the universal threshold: $\lambda=\sigma\sqrt{2\log(n)}$). In practice, this variance is unknown and needs to be estimated. One could estimate $\sigma^2$ in an iterative way, i.e. with a backfitting algorithm. We propose instead a direct method of estimation based on a QR decomposition of the linear part.

In wavelet approaches for standard nonparametric regression, a popular and well behaved estimator for the unknown standard deviation of the noise is the median absolute deviation (MAD) of 
the finest detail coefficients of the response divided by 0.6745 (see \shortciteA{Donoho96}). The use of the MAD makes sense provided that the wavelet representation of the signal to be denoised is sparse. However, such an estimation procedure cannot be applied without some pretreatment of the data in a partially linear model because the wavelet representation of the linear part of a PLM may be not sparse. Indeed, in practice we have observed that for many partly linear models such a procedure leads to biased estimations.

A QR decomposition on the regression matrix of the PLM allows to eliminate this bias. Since often the function wavelet coefficients at weak resolutions are not sparse, we only consider the wavelet representation at level $J=\log_2(n)$. Let $A_J$ be the wavelet representation of the design matrix $X$ at level $J$. The QR decomposition ensures that there exist an orthogonal matrix $Q$ and an upper triangular matrix $R$ such that $$A_J=Q\begin{pmatrix}R\\ 0\end{pmatrix}.$$ 

If $\bZ_J$, $\btheta_{0,J}$ and $\bepsilon_J$ denote respectively the vector of the wavelets coefficients at resolution $J$ of $Y$, $\bff$ and $U$, model (\ref{Wmodele}) gives $$Q^T \bz_J = \begin{pmatrix}R\\ 0\end{pmatrix}\bbeta_0+Q^T\btheta_{0,J}+Q^T\bepsilon_J.$$
It is easy to see that applying the MAD estimation on the last components of $Q^T \bz_J$ rather than on $\bz_J$ will lead to a satisfactory estimation of $\sigma$. Indeed thanks to the QR decomposition the linear part does not appear anymore in the estimation and thus the framework is similar to the one used in nonparametric regression. Following~\shortciteA{Donoho96}, the sparsity of the functional part representation ensures good properties of the resulting estimator.\\

\section{Simulation study}

The purpose of this section is to study through simulations several algorithms for estimating the linear part of a PLM model but also to evaluate the performance of the proposed estimators. Our wavelet estimation method for PLM will be also compared with a wavelet backfitting algorithm proposed by \shortciteA{FadiliBullmore}. 
As already noted, our estimation method allows us to first estimate the linear regression parameter vector $\bbeta_0$ independently of the nonparametric part, and to then proceed to the estimation of the functional part of the PLM model. The $M$-estimation $\hat\bbeta_n$ of $\bbeta_0$ is obtained by means of iterative optimization procedures that are more or less efficient, but usually much faster than backfitting procedures, as we shall see. Before proceeding to the analysis of our simulation results, we briefly recall two particular optimization algorithms that may be used for estimating the linear part.

\subsection{Half-quadratic algorithms}

The minimization problem we have to solve is of the form:
\begin{equation}\label{mini}
\hat{\bbeta}_n=\argmin_{\bbeta} J(\bbeta) ~~\text{with}~~J(\bbeta)=\sum_{i=1}^n\rho_\lambda(z_i-\bA_i^T\bbeta).
\end{equation}
Minimizers of $J(\bbeta)$ can be obtained using standard optimization tools such as relaxation, gradient, conjugated gradient and so on, but even if the loss function $\rho_\lambda$ is convex, its second derivative is large near to zero, so the optimization may be slow. For this reason, specialized optimization schemes have been conceived.
A very successful approach is {\em half-quadratic optimization}, proposed in \shortciteA{GemanReynolds} and \shortciteA{GemanYang} for cost functions of the above form. The idea is to associate with every $\bbeta$ in (\ref{mini})
an auxiliary variable ${\bf c}$ and to construct an augmented criterion $K$, such that for every $\mathbf{c}$ fixed, the function $\bbeta\to K(\bbeta,\mathbf{c})$ is quadratic (hence quadratic programming can be used) whereas for every $\bbeta$ fixed, each $\mathbf{c}$ can be computed independently using an explicit formula. The augmented criterion $K$ is chosen to have the same minimum as $J$, attained for the same value of $\bbeta$. The optimization problem of the augmented energy can be solved iteratively. 
At each iteration one realizes an optimization with respect to $\bbeta$ for $\mathbf{c}$ fixed and a second with respect to $\mathbf{c}$ for $\bbeta$ fixed. More precisely, if $\bbeta^{(m)}$ and ${\bf c}^{(m)}$ are the values given after $m$ iterations, the $(m+1)^{th}$ step of the algorithm actualizes these values through: \begin{equation}\label{iteration}\begin{array}{ccc}
{\bbeta}^{(m+1)}& = & \argmin_{\bbeta} K(\bbeta,{\bf c}^{(m)})\\
{\bf c}^{(m+1)} & = & \argmin_{\bf c} K(\bbeta^{(m+1)},{\bf c}) 
\end{array} \end{equation}

This procedure leads to two algorithms, namely ARTUR and LEGEND, that are also referenced in the literature as IRLS and IMR. We refer to \shortciteA{NikolovaNg} for some theory on the their use with Huber M-estimation. These algorithms are used for example in robust recognition (see e.g. \shortciteA{Charbonnier},~\shortciteA{Dahyot-bayes} or~\shortciteA{Vik}). \shortciteA{Vik} in particular stresses the link between ARTUR and LEGEND and Huber's approach.

\subsubsection*{ARTUR}

The algorithm described hereafter is referenced as the ARTUR algorithm in the optimization literature
or as \textit{Iterative Reweighted Least Squares} (IRLS) in the robustness literature. Geman and Reynolds's theorem leads to an augmented criterion of the form 
$$K(\bbeta,{\bf c})=\sum_{i=1}^n c_i (z_i-\bA_i^T\bbeta)^2 +\Psi({\bf c}).$$ 
The auxiliary variable ${\bf c}$ corresponds to a weight on the residuals of the least squares fit, thus explaining the IRLS terminology. Intuitively, weights on large residuals have a tendency to eliminate
the corresponding responses from the fit. For $\bbeta$ fixed, the minimum is reached for $c_i=\frac{\rho_\lambda'(r_i)}{r_i}$ where $r_i$ is the $i$th residual $r_i=z_i-\bA_i^T\bbeta$. At this point the value of $\Psi$ is $\rho_\lambda(r_i)-\rho_\lambda'(r_i)r_i/2$.

The $m+1$ step of the ARTUR algorithm can therefore be described as follows: 
$$\left\{\begin{array}{rcl}
r_i^{(m)} & = & z_i-\bA_i^T\bbeta^{(m)}\\
c_i^{(m+1)} & = & \frac{\rho_\lambda'(2r_i^{(m)})}{2r_i^{(m)}}, \qquad \qquad \forall i\in\{1,\ldots, n\}\\
\bbeta^{(m+1)} & = & (A^T{\bf c}^{(m+1)}A)^{-1}A^T{\bf c}^{(m+1)}\bZ\\
\end{array}\right. $$

\subsubsection*{LEGEND}

LEGEND, or \textit{Iterative Modified Residuals} (IMR), is a slightly different algorithm. The auxiliary variable doesn't weight the residuals anymore but subtracts the larger values of the residuals instead. The existence of the corresponding augmented energy functional follows from the second theorem of \shortciteA{GemanReynolds}. The criterion to be minimized can be written as 
 $$K(\bbeta,{\bf c})=\sum_{i=1}^n (z_i-\bA_i^T\bbeta-c_i)^2 +\xi({\bf c}).$$
For $\bbeta$ fixed, the minimum is reached for $c_i=r_i\left(1-\frac{\rho_\lambda'(r_i)}{2r_i}\right)$ where $r_i$ $i$th residual $r_i=z_i-\bA_i^T\bbeta$ and at this point the function $\xi$ takes the value $\rho_\lambda(r_i)-\rho_\lambda'(r_i)^2/4$.

With similar notation as for the ARTUR algorithm, the $m+1$ step of the LEGEND algorithm can be described as follows: 
$$\left\{\begin{array}{rcl}
r^{(m)} & = & \bZ-A\bbeta^{(m)}\\
c_i^{(m+1)} & = & r_i^{(m)}\left(1-\frac{\rho_\lambda'(2r_i^{(m)})}{2r_i^{(m)}}\right)\qquad \qquad \forall i\in\{1, \ldots, n\}\\
\bbeta^{(m+1)} & = & (A^TA)^{-1}A^T(\bZ-{\bf c}^{(m+1)})\\
\end{array}\right. $$

Both ARTUR and LEGEND are very easy to program. \shortciteA{NikolovaNg} show that the risk obtained via the multiplicative form ARTUR is always smaller than the one obtained via the additive form, but the later one is numerically faster. The main reason for this is that under the multiplicative form a matrix inversion is performed within each iteration.

\subsection{Numerical simulations}

In this subsection, we give some simulation results. All the calculations were carried out in MATLAB 7.0 on a unix environment. For the DWT, we used the WaveLab toolbox developed by Donoho and his collaborators at the Statistics Department of Stanford University (http://www- stat.stanford.edu/$\tilde{}$wavelab). For each of the simulated examples in the sequel, we may summarize the various ingredients of our fitting procedure as follows:

\begin{enumerate}
\item Application on the observed data of the discrete wavelet transform (DWT) using the pyramidal algorithm of \shortciteA{Mallat};
\item Estimation of the variance $\sigma^2$ by means of a QR decomposition on the matrix of wavelet coefficients at maximal resolution followed by a MAD estimation;
\item \label{artleg1} Estimation of $\bbeta_0$ with ARTUR or LEGEND, solving (\ref{beta});
\item \label{artleg2} Estimation of ${\btheta_0}$ by soft thresholding of $\bZ-A\hat\bbeta_n$, given by (\ref{theta});
\item Finally, estimation of $\hat f_n$ by applying the inverse DWT on $\hat\btheta_n$.
\end{enumerate}
We will compare with Fadili and Bullmore's procedure that estimates conjointly $\bbeta_0$ and $\btheta_0$ using a backfitting algorithm.\\
In order to reduce the number of iterations, we have used a stopping criterion in both ARTUR et LEGEND: while fixing a larger upper bound for the total number of iterations allowed, we also consider that the algorithm has converged as soon as the difference between two successive iterations is smaller than some given threshold $\delta$. More precisely, the iterations are stopped as soon as $\frac{\|\bbeta^{(m+1)}-\bbeta^{(m)}\|_2}{\|\bbeta^{(m)}\|_2}<\delta$ or whenever we attain their upper limit. 

For illustration, we generated three test problems as follows. The nonparametric component $f_0$ was selected among two different functions, one sinusoidal function and one piecewise constant function. The covariate is chosen as $\bX_i=g(i/n)+\eta_i$ with polynomial functions $g$ and with the $(\eta_i)_{i=1, \ldots, n}$ generated independently from a centered distribution with finite variance, as explained in Section \ref{hypotheses}. For DWT, the filter we used is the Daubechies Symmlet filter with 8 vanishing moments.
The sample size we took was $n=2^8$. For each setting, 500 replicates of data with different $X$ and $u$ were generated. The variance of the noise was chosen such as the signal-to-noise ratios of the nonparametric and parametric component respectively were equal to 2.2 and 4.38. Such choices seem reasonable. With the simulated data, we then used the proposed algorithms to estimate the unknown parameters. For wavelet thresholding the universal threshold was used, while the termination tolerance $\delta$ was set to $10^{-5}$ for ARTUR and $10^{-10}$ for LEGEND. For \textit{Backfitting}, we have used the algorithm of \shortciteA{FadiliBullmore} with a tolerance level $\delta$ equal to $10^{-20}$. To save computational time we have also specified an upper limit of 2000 for the maximum number of iterations allowed.

\subsubsection*{Example 1: Sinusoidal test function}

In examples 1 and 2, the covariate was generated using the polynomial function $g(t)=t^5+2t$ and with the $(\eta_i)_{i=1, \ldots, n}$ generated independently from $N (0, 1)$. We have also run some numerical simulations with different design functions $g$ such as $g(t)=2^t$, $g(t)=e^{-t^2}$ or $g(t)=cos(t)$ with similar results, not reported here by the lack of space. It seems that assumption (A4) is not really necessary for asymptotic consistency.

We first consider the case of a sinusoidal function for the nonparametric part. In such a case one could obviously use smoothing splines based semiparametric estimation but it is interesting to see how our wavelet based procedure behaves. Figure~\ref{regressionb} displays the wavelet transform of the data and of the design matrix. Note that the sparse representation of the nonparametric part allows an efficient reduction of the bias between the observations and a linear model. The dashed lines in the plot displayed in Figure \ref{regressionb}, represent the lines $X_i\bbeta_0\pm \lambda$. Observations lying far out from these lines do not affect the estimation of $\bbeta_0$. 

\figs{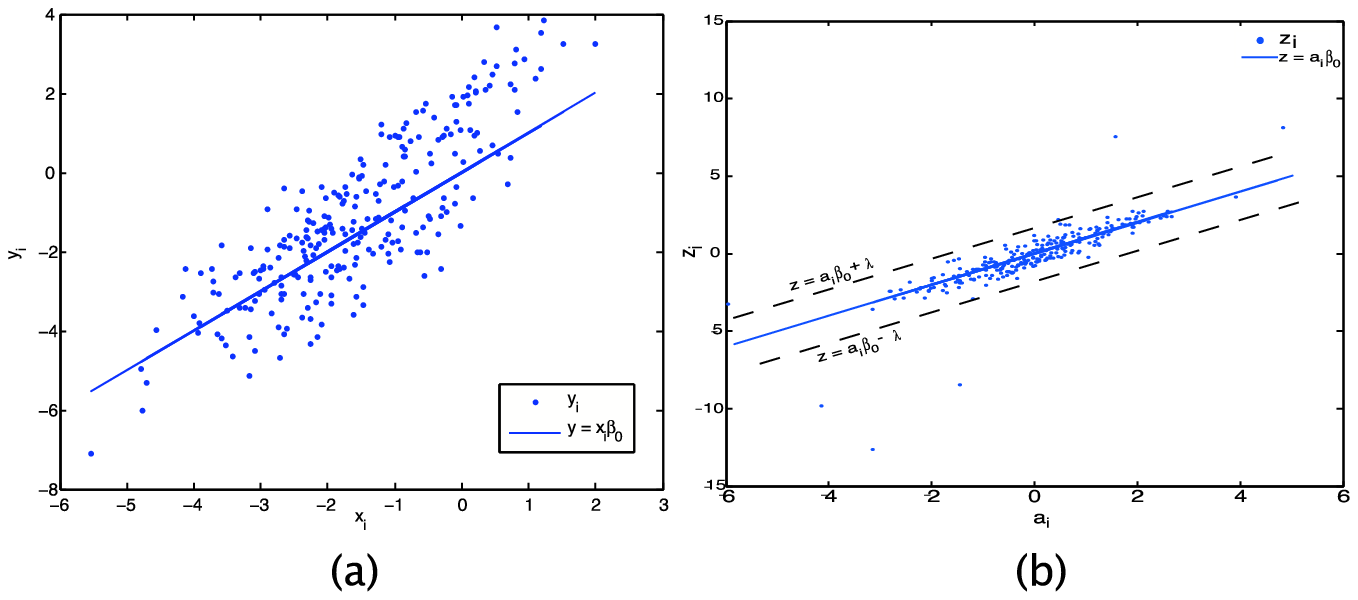}{12cm}{Wavelet transform of the data. Figure (a) represents the scatter plot of the observations $y_i$ versus the values of the covariates $X_i$. The line is the linear part of the model, of equation $y_i=X_i\bbeta_0$. Figure (b) is the scatter plot in (a) after the Discrete Wavelet Transform: it represents the coefficients $z_i$ versus $\bA_i$. The solid line is the linear part of the model (equation $z_i=A_i\bbeta_0$) and the dashed lines are the lines of equations $z_i=A_i\bbeta_0\pm \lambda$.}


We now evaluate the effect of the QR decomposition on the estimation of the noise, and we compare the computational time required by each of the algorithms, namely ARTUR, LEGEND and \textit{Backfitting} over the 500 replications of the experiment.

\vspace{0.5cm}

\begin{table}[htbp]
 \centering
 \begin{tabular}{@{} lcr @{}} 
 \toprule
 \multicolumn{3}{c}{Estimation of $\sigma$ by MAD} \\
 \cmidrule(r){1-3} 
 True value & without QR & with QR\\
 \midrule
 0.5 & 1.2222(0.0955) & 0.5023(0.0511) \\
 \bottomrule
 \end{tabular}
 \caption{{\small The mean values of the estimates and their standard deviation over the 500 simulations in Example 1 with $n=2^{8}$ (the standard deviation appears in brackets).}}
 \label{tab:table1}
\end{table}
\vspace{0.5cm}

From Table~\ref{tab:table1}, we get a fairly good impression on the effect of the QR decomposition on the estimation of the noise variance: 
the presence of the linear part introduces a strong bias in the MAD estimator, bias which is strongly diminished when using the QR decomposition. This also explains why in the comparison of their various thresholded estimators, \shortciteA{FadiliBullmore} often obtain estimators that are over-smoothed, since the variance that is used in their thresholds is over estimated. To be fair, we therefore
have adopted for all methods the universal threshold $\lambda=\sigma\sqrt{2\log(n)}$ with $\sigma$ estimated by MAD after a QR decomposition.

\vspace{0.5cm}

\begin{table}[htbp]
 \centering
 \begin{tabular}{@{} lccr @{}} 
  \toprule
  \multicolumn{4}{c}{Estimation of $\bbeta_0$} \\
  \cmidrule(r){1-4} 
  True value & \textit{Backfitting} & ~~ARTUR ~~& ~~LEGEND \\
  \midrule
  1 & 0.9000(0.0273) & 0.9417(0.0327) & 0.9417(0.0327) \\
  \midrule 
  Average computing time & 0.0936 & 0.0232 & 0.0151 \\ 
  \bottomrule
 \end{tabular}
 \caption{{\small The mean values of the estimates and their standard deviation over the 500 simulations in Example 1 (standard deviation appears in brackets) with $n=2^{8}$. The average MISE for the nonparametric part for these simulations is 0.1029 for ARTUR and LEGEND and 0.1098 for \textit{Backfitting}.}}
 \label{tab:table2}
\end{table}

From the last row of Table~\ref{tab:table2} one can see that both half-quadratic procedures (ARTUR and LEGEND) are faster than \textit{Backfitting} and the quality of estimation of both the parametric and nonparametric parts in terms on mean squared error is also better. The differences observed
in estimating $\bbeta_0$ between the various procedures is mainly due to the different tolerance levels
$\delta$ used by each. Note also that \textit{Backfitting} always stops because the maximal number of iterations is reached. The estimation given by \textit{Backfitting} could be improved but at the cost of a much larger computational time. 

\begin{figure}[!ht]
\begin{center}
\hspace{1.5cm} \includegraphics[width=12cm]{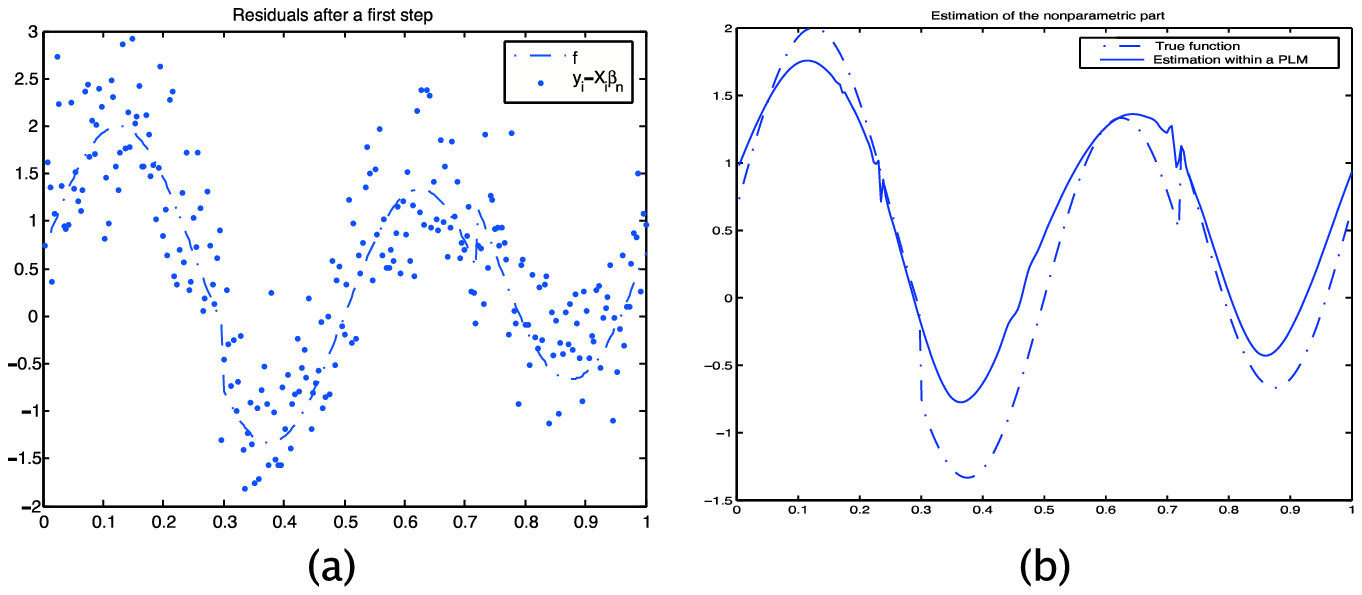}
\\
\end{center}
\vspace{-0.5cm}\caption{{\small Estimation of the nonparametric part in Example 1. Figure (a) represents the residuals obtained after estimation of the linear part of the models, meaning $z_i-\bA_i\hat\bbeta_n$, and the true functionnal part (dash). in Figure (b) we have the resulting estimation of the function (solid) and the true function (dash).}} \label{fct1}
\end{figure}

Recall that for both half-quadratic based algorithms, once the unknown parameter $\bbeta_0$ is estimated, a nonparametric wavelet based estimation procedure is applied to the resulting residuals
$y_i-\bX_i\hat\bbeta_n$ for estimation of the nonparametric part. Figure \ref{fct1} displays a typical example of these residuals and of the corresponding nonparametric estimation using ARTUR on one replication. 

For the value of the signal-to-noise ratio ($SNR_f=2.2$) adopted in our simulations for the nonparametric part, the estimator does not detect the discontinuity. However it produces results very similar to those by standard wavelet denoising of an identical nonparametric signal (without a linear part) with the same SNR, supporting our claim that the presence of the linear part in a PLM doesn't affect the estimation of the nonparametric part. 

\figs{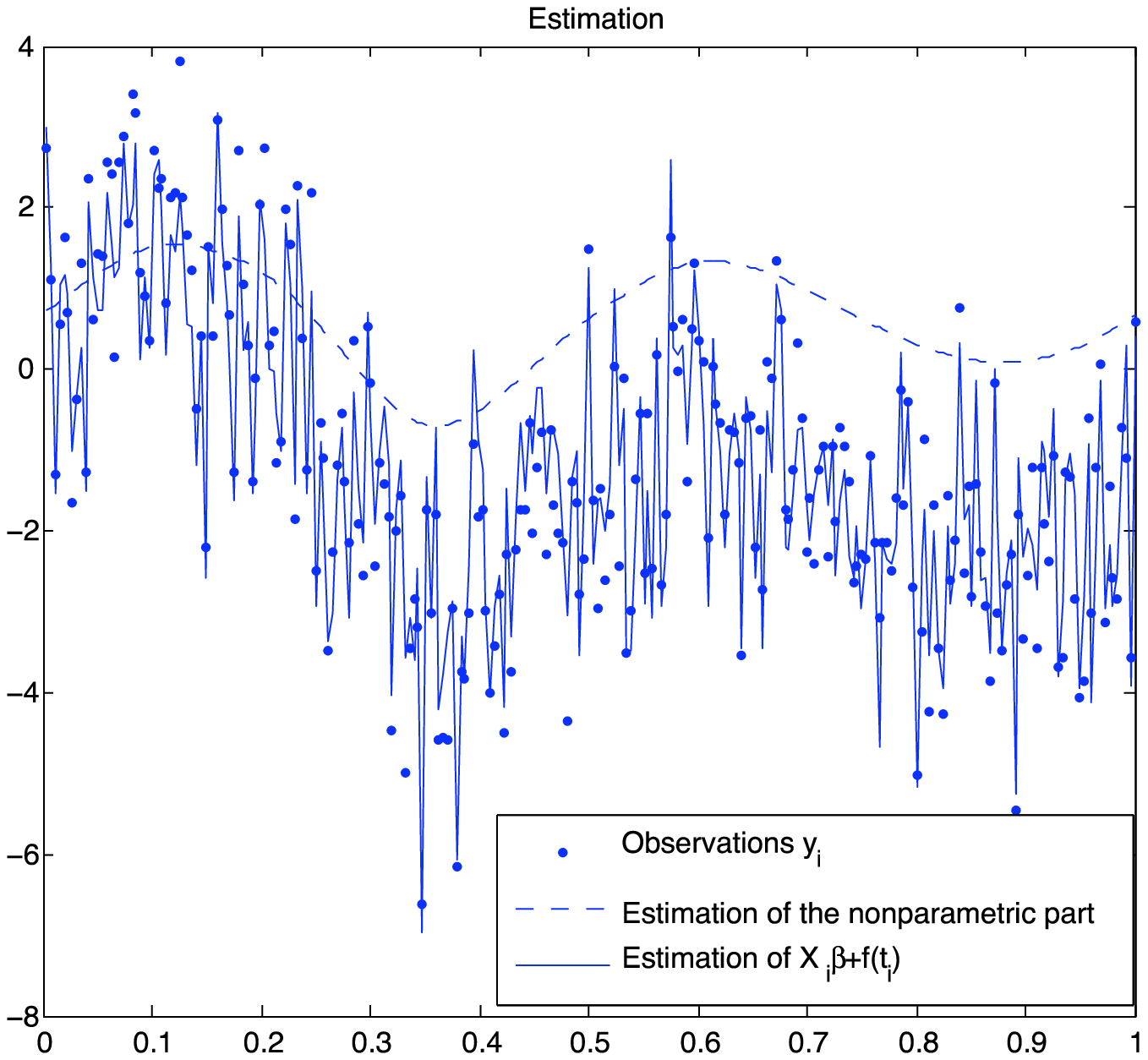}{6cm}{{\small A typical partial linear fit from Example 1. The figure represents the scatter plot of the observations, the estimated functionnal part (dash) and the parial linear fit (solid) for one of the simulation.}}

In their numerical implementation of ARTUR et LEGEND, both \shortciteA{Vik} and \shortciteA{Dahyot-comp} conclude that LEGEND converges faster, supporting the theoretical results of \shortciteA{NikolovaNg}. To share some light on this fact we have run some simulations with a larger number sample size. With $n=2^{10}$ observations and the same signal-to-noise ratio as before one can see a clear difference in computational time among the two algorithm for estimators with equivalent qualities,
as reported in Table~\ref{tab:table3}. 

\vspace{0.5cm}

\begin{table}[htbp]
  \centering
  \begin{tabular}{@{} lcr @{}} 
   \toprule
   \multicolumn{3}{c}{Estimation of $\bbeta_0$} \\
   \cmidrule(r){1-3} 
   True value & ~~ARTUR ~~ & ~~ LEGEND \\
   \midrule
    1 & 0.9762(0.0127) & 0.9762(0.0127) \\
    \midrule
    Average computing time & 0.2331 & 0.0166 \\
    Average number of iterations & 7 & 59 \\
   \bottomrule
  \end{tabular}
  \caption{{\small The mean values of the estimates and their standard deviation over the 500 simulations in Example 1 (the standard deviation appears in brackets) with $n=2^{10}$. LEGEND is much faster
  than ARTUR.}}
  \label{tab:table3}
\end{table}

\subsubsection*{Example 2: piecewise linear function}

We would like now to illustrate our estimation procedure when the nonparametric part is highly non regular. We thus consider a function $f_0$ which is piecewise constant. It is obvious that for such a function, our wavelet based procedure is better suited than a spline based procedure. All other setting adopted for these simulations are the same as those for example 1.

\begin{table}[htbp]
 \centering
 \begin{tabular}{@{} lcr @{}} 
 \toprule
 \multicolumn{2}{c}{Estimation of $\sigma$ by MAD} \\
 \cmidrule(r){1-2} 
 True value &  with QR\\
 \midrule
 0.5 &  0.49961(0.052741) \\
 \bottomrule
 \end{tabular}
 \caption{{\small The mean values of the estimates and their standard deviation over the 500 simulations in Example 2 (the standard deviation appears in brackets) for $n=2^8$.}}
 \label{tab:table4}
\end{table}

\vspace{0.5cm}

The results given in Table \ref{tab:table5} reinforce our claim from example 1 that half-quadratic algorithms are more efficient than {\it Backfitting}. Note moreover that the non regularity of the nonparametric part does not seem to affect the quality of the estimation of the vector of regression parameters.

\vspace{0.5cm}

\begin{table}[htbp]
 \centering
 \begin{tabular}{@{} lccr @{}} 
  \toprule
  \multicolumn{4}{c}{Estimation of $\bbeta_0$ for $n=2^8$} \\
  \cmidrule(r){1-4} 
  True value & \textit{Backfitting} & ~~ARTUR ~~ & ~~ LEGEND \\
  \midrule
  1 & 0.8999(0.0273) & 0.9548(0.0309) & 0.9548(0.0309) \\
  \midrule 
  Average computing time & 0.0744 & 0.0209 & 0.0139 \\ 
  \bottomrule
 \end{tabular}
 \caption{{\small The mean values of the estimates and their standard deviation over the 500 simulations in Example 2 (standard deviation appears in brackets). The average MISE for the nonparametric part for these simulations is 0.1012  for ARTUR and LEGEND and 0.1078 for \textit{Backfitting}.}}
 \label{tab:table5}
\end{table}

\begin{table}[htbp]
 \centering
 \begin{tabular}{@{} lccr @{}} 
  \toprule
  \multicolumn{3}{c}{Estimation of $\bbeta_0$ for $n=2^{10}$} \\
  \cmidrule(r){1-3} 
  True value & ~~ARTUR ~~ & ~~ LEGEND \\
  \midrule
  1 & 0.9554(0.0149) &  0.9554(0.0149) \\
  \midrule 
  Average computing time & 0.3036 & 0.0209 \\ 
  \bottomrule
 \end{tabular}
 \caption{{\small The mean values of the estimates and their standard deviation over the 500 simulations in Example 2 (standard deviation appears in brackets). The average MISE for the nonparametric part for these simulations is  0.0584 for ARTUR and LEGEND.}}
 \label{tab:table6}
\end{table}

As in example 1, one can see from Table \ref{tab:table5} and Table \ref{tab:table6} that LEGEND outperforms ARTUR, and that the difference of computing time increases with the number of observations $n$.

\begin{figure}[!ht]
\begin{center}
\hspace{1.5cm} \includegraphics[width=12cm]{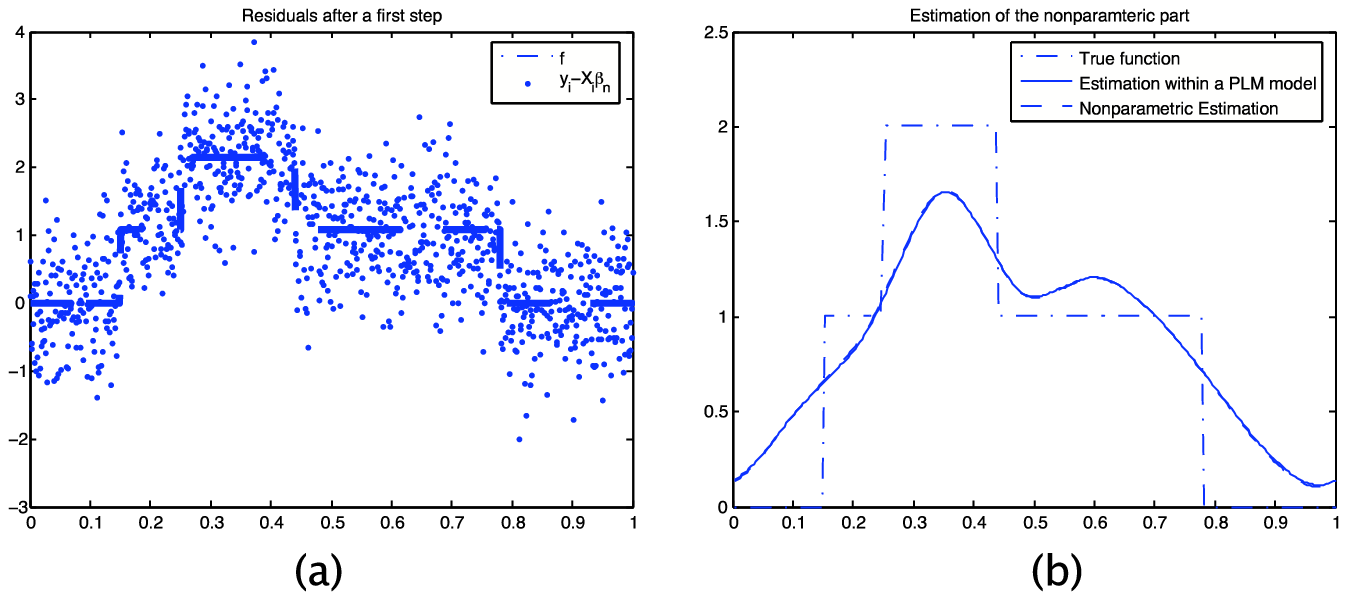}
\\
\end{center}
\vspace{-0.5cm}\caption{{\small Estimation of the nonparametric part in Example 2. Figure (a) represents the residuals obtained after estimation of the linear part of the models, meaning $z_i-\bA_i\hat\bbeta_n$, and the true functionnal part (dash). In Figure (b) we have the resulting estimation of the function (solid) and the true function (dash).}} \label{fct2}
\end{figure}

The estimation of the nonparametric part does not detect the discontinuities of the function. Yet compared to standard wavelet denoising in a nonparametric regression model with the same SNR, the estimation obtained in the PLM is very similar. The bad visual quality of the estimation results from the choice of the signal-to-noise ratio ($SNR_f=2.2$) adopted in our simulations rather than the presence of the linear part.

\subsubsection*{Example 3: dimension 4}

We now consider a case where the vector of parameter $\bbeta$ belongs to $\R^4$ (the dimension of the design regression matrix $X$ is then $n\times 4$). The nonparametric part $f_0$ is the same as in example 2, meaning that the function is highly irregular. The SNR for the global model was chosen equal to 5.99, with a SNR equal to 4.38 for the nonlinear part. One may summarize the results for this example in the above tables. 

\begin{table}[htbp]
  \centering
  \begin{tabular}{@{} lr @{}} 
   \toprule
   \multicolumn{2}{c}{Estimation of $\sigma$ by MAD with QR} \\
   \cmidrule(r){1-2} 
   True value & with QR\\
   \midrule
    0.5 & 0.52261(0.053808) \\ 
   \bottomrule
  \end{tabular}
  \caption{{\small The mean values of the estimates and their standard deviation over the 500 simulations in Example 3 (the standard deviation appears in brackets).}}
  \label{tab:table7}
\end{table}
 
\vspace{0.5cm}

\begin{table}[htbp]
  \centering
  \begin{tabular}{@{} lccr @{}} 
   \toprule
   \multicolumn{4}{c}{Estimation of $\bbeta_0$} \\
   \cmidrule(r){1-4} 
   True value & \textit{Backfitting} &~~ARTUR ~~ & ~~ LEGEND \\
   \midrule
	-1  &   -1.4969(0.45822) & -0.7203(0.461)   & -0.7203(0.461)  \\
	3   &    2.8563(0.09770) &  2.9168(0.09941) &  2.9168(0.09941) \\ 
	0   &   -0.1201(0.33685) &  0.0125(0.34415) &  0.0125(0.34415) \\  
	8   &    7.5601(0.16772) &  7.7112(0.18525) &  7.7112(0.18525)  \\
    \midrule
    Mean squared error & 0.8434  &  0.5438  &  0.5438 \\ 
     \midrule
    Average computing time & 0.1602 & 0.0305 & 0.0234 \\
   \bottomrule
  \end{tabular}
  \caption{{\small The mean values of the estimates and their standard deviation over the 500 simulations in Example 3 (the standard deviation appears in brackets) for a given value of the true $\bbeta_0$. The average MISE for the nonparametric part for these simulations is 0.2140 for ARTUR and LEGEND and 0.2164 for \textit{Backfitting}.}}
  \label{tab:table8}
\end{table}

As one can see with computational times that are similar for all procedures, both half-quadratic algorithms outperform \textit{Backfitting} in terms of the MSE.

\begin{figure}[!ht]
\begin{center}
\hspace{1.5cm} \includegraphics[width=12cm]{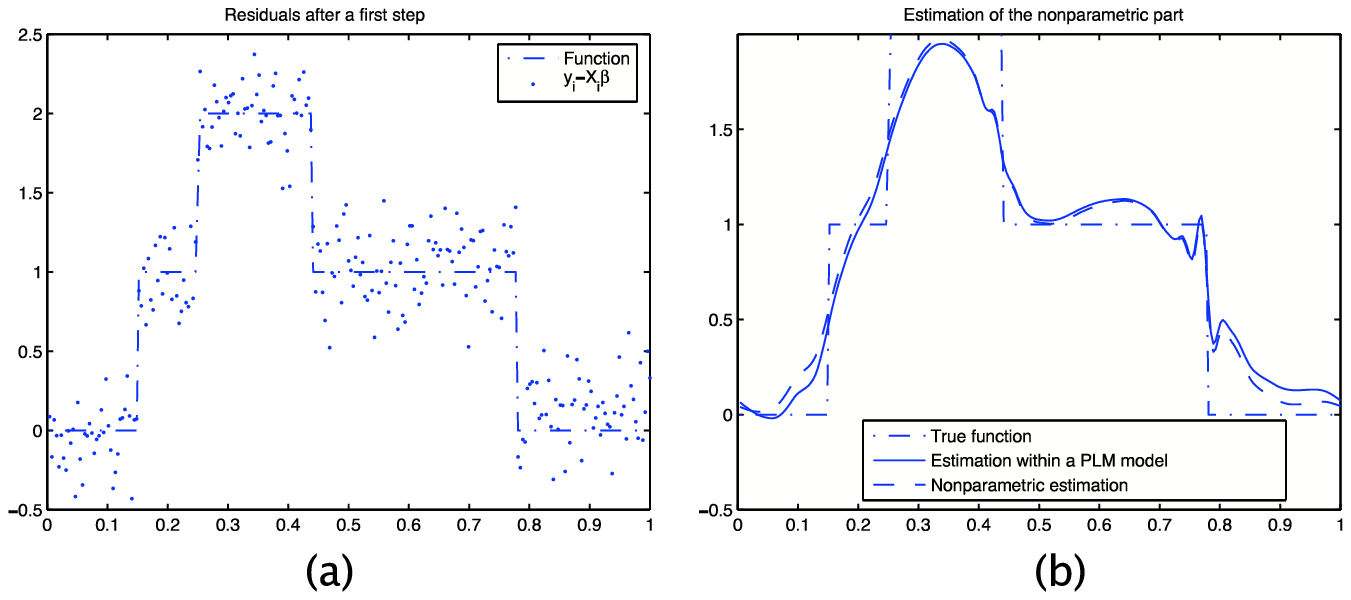}
\\
\end{center}
\vspace{-0.5cm}\caption{{\small Estimation of the nonparametric part in Example 3. Figure (a) represents the residuals obtained after estimation of the linear part of the models, meaning $z_i-\bA_i\hat\bbeta_n$, and the true functionnal part (dash). In Figure (b) we have the resulting estimation of the function (solid) and the true function (dash).}} 
\end{figure}

\begin{figure}[!h]
\begin{center}
\hspace{1.5cm} \includegraphics[angle=0,scale=0.45,height=5cm,width=6cm]{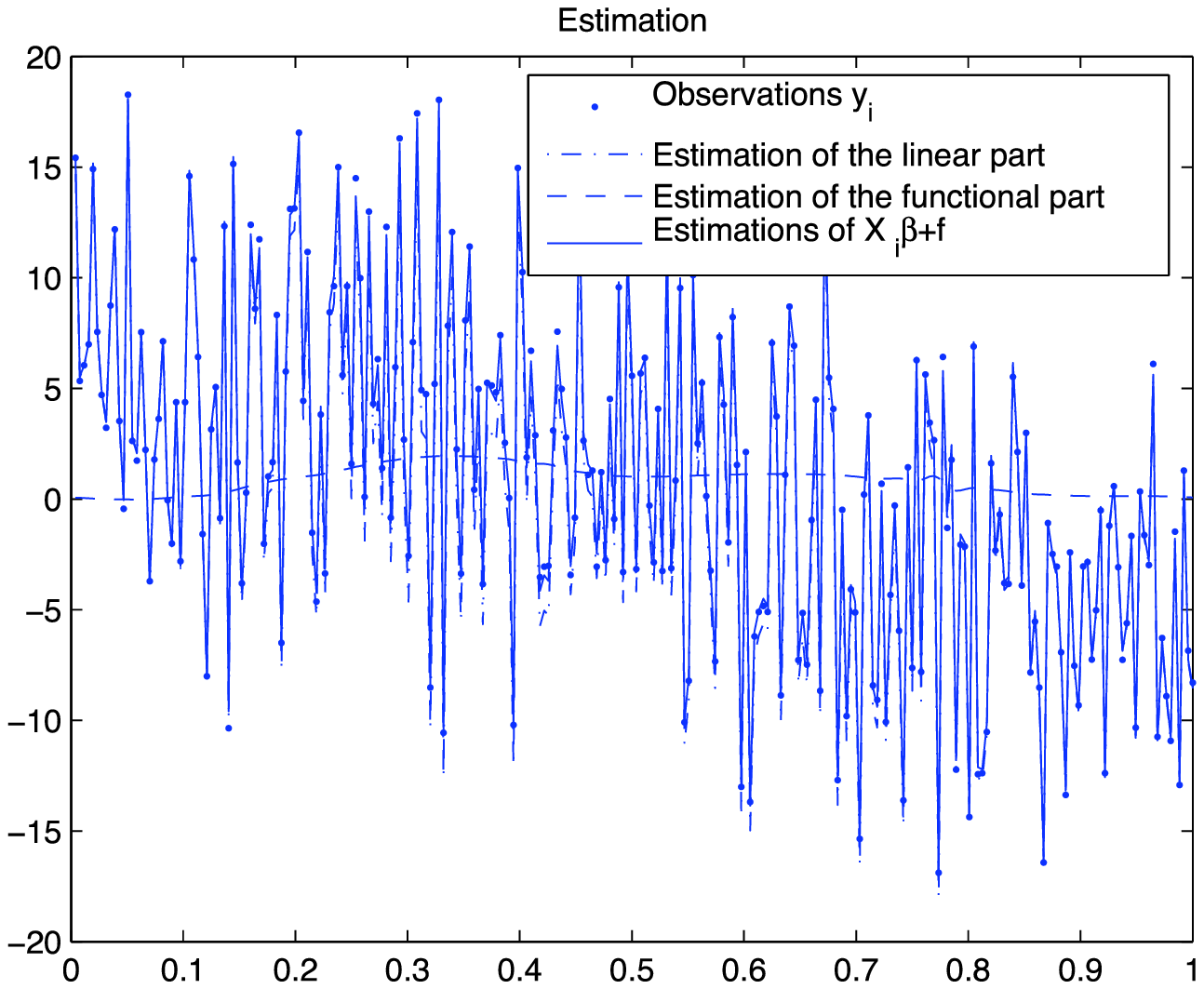}
\\
\vspace{-0.5cm}\caption{{\small A typical partial linear fit from Example 3. The figure represents the scatter plot of the observations, the estimated functionnal part (dash) and the partial linear fit (solid) for one of the simulation.}}
\end{center}
\end{figure}

As for examples 1 and 2, when the sample size increases, among the half-quadratic algorithms the LEGEND one is much faster.


\section*{Conclusion}

This paper develops a powerful penalized least squares estimation in partially linear models, based on a wavelet expansion of the nonparametric part. Choosing an appropriate penalty on the wavelet coefficients of the function, the procedure leads to an estimation of the linear part of partly linear models independent from the nonparametric part, while the estimation of the nonparametric part is adaptative relatively to the smoothness of the function. Since the functionnal part of the model has a sparse representation, the estimation of the regression parameters vector is moreover interpreted as a common M-estimation. In the particular case of an $l^1$-penalty (leading to soft thresholding and Huber's estimator) the near-minimaxity of the estimation of both parametric and nonparametric parts of a partially linear model is established, and the result is avalaible for a large class of functions, including nonsmooth irregular functions. From an implementation point of view, half-quadratic algorithms are proposed that appear to give good results on simulation studies.

Our ongoing research is focusing on exploring the asymptotic properties of the procedure for other thresholding schemes and in more general frameworks such as nonequidistant designs for the nonparametric part.

\section*{Acknowledgements}
Part of this work was supported by the `IAP Research Network
P5/24'. The author would like to
thank Dr. Fadili and Dr. Bullmore for kindly providing the Matlab codes implementing the backfitting procedures used in the paper.

\subsection*{}

\section{Appendix}
\subsection*{Appendix A. Discussion of the assumptions.}\label{hypotheses}

In this Section, we study wether the assumptions made in Theorem \ref{resultat} are reasonable in practice. Following \shortciteA{Rice} or \shortciteA{Speckman} we suppose that the design matrix $X$ can be written as a sum of a deterministic function and a noise term. The $(i,j)$-component of $X$ can be written as $x_{i,j}=g_i(t_j)+\xi_{i,j}$ with functions $g_i$ such that $\int f g_i=0$ and where $\xi_{i,j}$ denotes a realization of a random variable $\xi_i$. The variables $(\xi_i)_{i=1, \ldots, n}$ are supposed to be independent and identically distributed, centered and with finite variance, independent from the $u_i$. With these notation, assumptions (A1), (A2) and (A4) become:
\begin{description} 
\item[(A1)] The norm of $\frac 1 n X^T\bff$ can be decomposed as follows:\\ $\|\frac 1 n X^T\bff_0\|^2=\sum_{j=1}^p\left(\frac1 n \sum_{j=1}^n g_{i}(t_j)f(t_j)+\frac1 n \sum_{j=1}^n \xi_{i,j}f(t_j)\right)^2. $ \\ The convergence towards 0 of the first term is ensured by the assumption that $\int f g_i=0$ for all $i=1, \ldots, p$. We can prove that the second term tends to 0 almost surely. 
\begin{rmq} When we suppose that $\forall i,~\int f g_i=0$, this impose that either the integral of $f$ is equal to zero or the vector $\1_{n\times 1}$ is not in the space spanned by the columns of $X$. This is the usual assumption for identifiability in PLM (e.g.~\shortciteA{Chen88} or~\shortciteA{DonaldNewey}).
\end{rmq}

\item[(A2)] Let $V(g)$ be the matrix with entries $\int g_i g_j$ and $V$ denotes the covariance matrix of the variables $(\xi_i)_{i=1, \ldots, n}$. One can prove that $\frac 1 n X^T X$ converges almost surely to $V(g)+V$. It is sufficient to assume that the family $(g_i)_{i=1, \ldots, n}$ is $\L^2$-orthogonal in order that the matrix $V(g)+V$ is non singular.

\item[(A4)] Actually, it is equivalent to prove that $\frac 1 n \sup\|\bA_i\|^2\to 0$ to get (A4). 
For $i\in\{1,\ldots, n\}$ given, $\frac 1 n \|\bA_i\|^2$ is equal to $n^{-1}\bA_i^T\bA_i=\sum_{l=1}^p\left[\frac{1}{{n}}\sum_{j=1}^n\psi_i(t_j)X_{j,l}\right]^2$. With the previous notation, $x_{j,l}=g_l(t_j)+\xi_{j,l}$ and we can establish that $n^{-1}\bA_i^T\bA_i$ tends almost surely to $\sum_{l=1}^p \left(\int \psi_i g_l\right)^2$. This can also be written as $n^{-1}\bA_i^T\bA_i\sim \sum_{l=1}^p (w^l_i)^2$ with $(w_i^l)_{i=1,\ldots, n}$ wavelets coefficients of the functions $g_l$. \\
If, for all $l=1, \ldots, p$, $g_l$ is a polynomial function whose degree is less than or equal to the number of vanishing moments $N$ of the wavelet mother, then this assumption holds.

\end{description} 

Hypothesis (A3) is not detailled here because even if it does not seem very constraining, it is difficult to study its feasibility.

To conclude, when the design $X_i$, $i=1, \ldots, n$ can be written as $X_i=g_i+\xi_i$ with $g_i$ orthogonal polynomial functions with a degree less than or equal to $N$, and with $\xi_i$ centered independent random variables with finite variance, whenever $\int f g_i=0$ for all $i$, assumptions (A1), (A2) and (A4) hold.

\subsection*{Appendix B. Proofs of the main results}


\subsubsection*{B.1. Preliminary result}

\begin{Proposition}\label{controle}~\\
When assumptions (A2) and (A3) hold, 
$$\frac{1}{\sqrt n} \somme_{i=i_0}^n \rho'_{\lambda}({\theta_0}_i+\epsilon_i)\bA_i=\bigcirc_\P(\lambda)$$

~
\end{Proposition}
This result comes from Bernstein's inequality applied to the random variables  $Y_{i,j}=\frac{A_{i,j}}{\sqrt{n}}\frac{\rho'_\lambda({\theta_0}_i+\epsilon_i)}{\lambda}$, $i=1, \ldots, n$, for any fixed $j$ in $\{1, \ldots, p\}$. Indeed, these variables are almost surely uniformly bounded and $\sum_{i=1}^n \E[Y_{i,j}^2]$ is bounded, due to the following lemma:

 \begin{Lemme}\label{Ai} If (A2) and (A3) hold, \begin{description} \item[(i)] $n^{-1/2}\sup_{i=1, \ldots, n} \|\bA_i\|\to 0$ \item[(ii)] $n^{-1}\sum_{i=1, \ldots, n} \|\bA_i\|^2=\bigcirc(1)$ \end{description} \end{Lemme} This result lies on the observation that $\|\bA_i\|^2=\bA_i^T(A^TA)^{1/2}(A^TA)^{-1/2}\bA_i,$ and consequently $\|\bA_i\|\leq n^{1/2} \|(\frac 1 n A^TA)^{1/2} \| h_i^{1/2}$.
 

\subsubsection*{B.2. Variables transform}

Let us recall that we are studying the model \begin{equation}\label{modele}
z_i=\bA_i^T\bbeta_0+\theta_{0i}+\epsilon_i \quad \text{under~(A2)-(A4)}
\end{equation}
(Assumption (A1) is an identifiability assumption and does not intervene in the proofs).
Following \shortciteA{Huber} or \shortciteA{BaiRaoWu}, we build an equivalent model by a change of variables. Let us define the following transforms:
\begin{eqnarray*}
R&=&A(A^TA)^{-1/2},\\
{\balpha}&=&\frac{1}{\lambda}(A^TA)^{1/2}({\bbeta}-{\bbeta}_0)\\
d_i &=&\frac{1}{\lambda}(\theta_i+\epsilon_i).
\end{eqnarray*}
The results may be established equivalently for the following model:
\begin{equation}\label{modelebis}
z_i=\bR_i^T\balpha_0+d_i \quad \text{under~(A2'')-(A4'')};
\end{equation}
\begin{description}
\item[(A2'')] $R^TR=I_p$.

\item[(A3'')] $h=\displaystyle{\max_{i=i_0, \ldots, n}} \bR_i^T\bR_i$ tends to 0.

\item[(A4'')] $K_n'':=\sum_{i=i_0}^n \bR_i\bR_i^T \E\left[\rho''_{1}(d_i)\right]$ tends to $K_0''$, non singular matrix.

\end{description}
As the Huber cost function has scale transform properties: \begin{equation}\label{homothetie} \text{for~any~} v>0,~~\rho_{\lambda}(u)=v^2\rho_{\lambda/v}(u/v),\end{equation} we then can prove that in the model (\ref{modelebis}), the estimator $\hat\alpha_n$ is solution of the minimization problem $$\hat\balpha_n=\argmin_{\balpha} \sum_{i=1}^n \rho_{1}(d_i-\bR_i^T\balpha).$$  
As $\rho'_{\lambda}(u)=\lambda\rho'_{1}(u/\lambda)$ and $\rho''_{\lambda}(u)=\rho''_{1}(u/\lambda)$, we have $K_0''\sim \Sigma^{-1}K_0$ and Proposition \ref{controle} becomes in (\ref{modelebis}): $$\sum_{i=i_0}^n \rho'_{1}(d_i)\bR_i=\bigcirc_\P(1).$$
In all the proofs, we will consider the model (\ref{modelebis}) and obtain the consistency results thanks to the mentionned transforms. 

\subsubsection*{B.3. Convergence of the criterion}

\begin{Proposition}\label{convergencecritere}~\\
Let $c$ be a strictly positive constant. Suppose (A1) to (A4) hold. Then, \begin{eqnarray*}
\sup_{\left\{\|\bbeta-\bbeta_0\|\leq c \lambda n^{-1/2}\right\}} & \hspace{-0.05cm}
\frac{1}{\lambda^2}\left|\sum_{i=i_0}^n \left(\rho_{\lambda}({\theta_0}_i+\epsilon_i-\bA_i^T(\bbeta-\bbeta_0))-\rho_{\lambda}({\theta_0}_i+\epsilon_i)\right)\right.\\
& \left. +\somme_{i=i_0}^n \rho'_{\lambda}({\theta_0}_i+\epsilon_i)\bA_i^T(\bbeta-\bbeta_0) - n \frac 1 2 (\bbeta-\bbeta_0)^T K_0(\bbeta-\bbeta_0)\right|\tend^\P 0.
\end{eqnarray*}
\end{Proposition}
The proof is built on two phases: we first approximate the Huber cost function $\rho$ with a smoother function, keeping a control on the third derivative; secondly, we develop a scheme of proof very similar to \shortciteA{BaiRaoWu} in the transformed model (\ref{modelebis}). The main argument is the convexity of $\rho$, which allows in particular the use of Rockafellar's theorems. 

\subsubsection*{B.3.1. Approximation of Huber cost function}

The approximation is built by three successive integrations.
Let $0<\delta<1$. We define $r^3_\delta$ on $\R$: $$r^3_\delta:u \mapsto \begin{cases} 
\frac{6}{\delta^3}(u-(1-\delta/2))(u-(1+\delta/2)) &\text{if~} 1-\delta/2<|u|<1+\delta/2\\ 
0 &\text{otherwise}\end{cases}.$$
We introduce next, $r^2_\delta$ primitive of $r^3_\delta$ equal to zero at $1+\delta/2$, $r^1_\delta$ primitive of $r^2_\delta$ equal to zero at $0$ and $r_\delta$, primitive of $r^1_\delta$ equal to zero at 0.

The function series $\tilde \rho_{1}=r_{1/n^2}^3$ is a series of convex functions $\mathcal C^3$, which converges uniformly towards $\rho_1$ when $n$ goes to infinity. We can furthermore prove that $\int|\tilde\rho^{(3)}_1| \leq 12$, and that
\begin{eqnarray} \label{rho0}
n\|\tilde\rho_1-\rho_1\|_\infty & \tend_{n\to\infty} & 0,\\ \label{rho1}
n \|\tilde\rho'_1-\rho'_1\|_\infty & \tend_{n\to\infty} & 0,\\
\|\tilde\rho''_1-\rho''_1\|_\infty & \leq & 1.
\end{eqnarray}
Moreover, $\tilde \rho''_1$ and $\rho''_1$ only differ from each others on two intervals of length $1/n^2$.

\subsubsection*{B.3.2. Preliminary tools}

\begin{Proposition}\label{rock1}
Let $C$ be an open compact set of $\R^m$. We consider $(f_n)_{n\in\N}$ and $f$ a family of convex functions defined on $C$ and taking their values in a given probability space $(\Omega,P,\mu)$. Suppose for all $u\in C$, $f_n(u)-f(u)$ converges in probability to $0$. Then the convergence in probability of $\sup_{\{u\in C\}}f_n(u)-f(u)$ towards $0$ is acquired.
\end{Proposition}

\begin{proof}
We recall a theorem given in \shortciteA{Rockafellar} (Theorem 10.8, page 90):
\begin{Proposition}
Let $\cal C$ be an open compact set of $\R^m$. We consider $(f_n)_{n\in\N}$ and $f$ a family of finite convex functions defined on $C$. Suppose the series $(f_n)_{n\in\N}$ converges simply to $f$ on $\cal C$. Then the convergence is uniform on $\cal C$.
\end{Proposition}

In order to obtain a similar result for the convergence in probability, we may use the following characterization of such a convergence: \begin{Lemme}\label{ps}
Let $(X_n)_n$ be a series of random variables and $X$ a random variable. The series $(X_n)$ converges in probability towards $X$ if and only if from all subsequence of $X_n$ we can extract a series which tends almost surely to $X$.
\end{Lemme}

Consider $f_{\nu(n)}$ a subsequence of $f_n$. We would like to find $\bbeta(n)$, subsequence of $\nu(n)$, such that for all $u\in \cal C$, $f_{\bbeta(n)}(u)-f(u) \tend^{a.s.} 0$. The Lemma \ref{ps} tells us that for all $u\in \cal C$ there exists $\eta_u(n)$ extraction of $\nu(n)$ such that $f_{\eta_u(n)}(u)-f(u) \tend^{a.s.} 0$. Let us consider $\mathcal D=\{u_0,u_1,u_2\ldots\}$ dense and countable subset of $\cal C$. Using a diagonal procedure, we can exhibit $(\bbeta(n))$ such that for all $u\in \mathcal D$, we have $f_{\bbeta(n)}(u)-f(u) \tend^{a.s.} 0$. Afterwards, the convergence of $f_{\bbeta(n)}-f$ on $\cal C$ holds by density of $\mathcal D$ and continuity of $f_{\bbeta(n)}-f$.
Applying Rockafellar's Theorem, we obtain that $\displaystyle{\sup_{u\in C}} f_{\bbeta(n)}(u)-f(u)$ tends almost surely to $0$.

To conclude, we have proved that from all subsequence $\displaystyle{\sup_{u\in C}} f_{\nu(n)}(u)-f(u)$ of $\displaystyle{\sup_{u\in C}} f_{n}(u)-f(u)$
we could extract a series which converges almost surely to 0. This finishes the proof using Lemma \ref{ps}. 
\end{proof}

\subsubsection*{B.3.3. Convergence criterion}

Let $c>0$. We are going to prove that in model (\ref{modelebis}) we have: \begin{equation}\label{result1}
\sup_{\left\{\|\balpha\|\leq c \right\}} 
\left|\sum_{i=i_0}^n \left(\rho_{1}(d_i-\bR_i^T\balpha)-\rho_{1}(e_i)\right)+\sum_{i=1}^n\rho'_1(d_i)\bR_i^T\balpha
- \frac 1 2 \balpha^T K_0''\balpha\right|\tend^\P 0.
\end{equation} Note that in the initial model (\ref{modele}), this is equivalent to \begin{eqnarray*}
\sup_{\left\{\|\bbeta-\bbeta_0\|\leq c \lambda n^{-1/2}\right\}} &
\frac{1}{\lambda^2}\left|\sum_{i=i_0}^n \left(\rho_{\lambda}({\theta_0}_i+\epsilon_i-\bA_i^T(\bbeta-\bbeta_0))-\rho_{\lambda}({\theta_0}_i+\epsilon_i)\right)\right.\\
& \left. +\somme_{i=i_0}^n \rho'_{\lambda}({\theta_0}_i+\epsilon_i)\bA_i^T(\bbeta-\bbeta_0) - n \frac 1 2 (\bbeta-\bbeta_0)^T K_0(\bbeta-\bbeta_0)\right|\tend^\P 0.
\end{eqnarray*}

\hspace{-0.5cm} $\bullet$ We introduce:
$$\Delta(\balpha):=\sum_{i=1}^n \left(\tilde\rho_{1}(d_i-\bR_i^T\balpha)-\tilde\rho_{1}(d_i)+\tilde\rho_{1}'(d_i)\bR_i^T\balpha\right).$$
The cost function $\tilde\rho_{1}$ is convex. For every $i$, it gives the upper bound:
\begin{equation}\left|\tilde\rho_{1}(d_i-\bR_i^T\balpha)-\tilde\rho_{1}(d_i)+\tilde\rho'_1(d_i)\bR_i^T\balpha\right| \label{ineg} \leq |\tilde\rho'_{1}(d_i-\bR_i^T\balpha)-\tilde\rho'_{1}(d_i)| |\bR_i^T\balpha| .\end{equation}
This inequality gives a bound of the variance of $\Delta(\balpha)$:
$$Var(\Delta(\balpha))\leq \sum_{i=1}^n \E\left[\left(\tilde\rho_{1}'(d_i-\bR_i^T\balpha)-\tilde\rho'_{1}(d_i)\right)^2\right] |\bR_i^T\balpha|^2.$$

The function $\tilde\rho'_{1}$ being 1-Lipschitz, $$
\forall n\in\N,\;\forall i=1, \ldots, n,\;\forall u\in\R^+,\;\E\left(\tilde\rho'_{1}(d_i+u)-\tilde\rho'_{1}(d_i)\right)^2\leq u^2.$$ Consequently, $$Var(\Delta(\balpha))\leq \sum_{i=1}^n|\bR_i^T\balpha|^4\leq \|\balpha\|^4\sum_{i=1}^n\|\bR_i\|^4.$$ As $\balpha$ is supposed to be bounded and $\sum_{i=1}^n |\bR_i|^4\leq h\sum h_i=hp$ tends to 0, we obtain that $Var(\Delta(\balpha))$ tends to 0.
Bienaym\'e-Tchebychev inequality ensures then that $|\Delta(\balpha)-\E\Delta(\balpha)|$ converges towards 0 in probability.

\smallskip
\hspace{-0.5cm} $\bullet$ {\bf The term $\E\Delta(\balpha)$.}

As the function $\tilde\rho$ is $\mathcal C^3$, the Taylor expansion of degree 2 with a rest of an integral form of $\tilde\rho_{1}$ on a neighborhood of $d_i$ exists. It gives:
\begin{eqnarray*}\tilde\rho_{1}(d_i-\bR_i^T\balpha)-\tilde\rho_{1}(d_i)+\tilde\rho'_1(d_i)\bR_i^T\balpha -\frac 1 2 \tilde\rho_{1}''(d_i)\balpha^T\bR_i\bR_i^T\balpha\\ = -\int{\tilde\rho_{1}^{(3)}(t)(d_i-t)^3\1_{d_i-\bR_i^T\balpha\leq t\leq d_i}dt/6}.\end{eqnarray*} 

Using the bound $\int\left|\tilde\rho_{1}^{(3)}(t)\right|dt\leq 12,$ obtained when constructing $\tilde\rho$, we obtain:
\begin{equation*}\E\left|\sum_{i=1}^n\left(\tilde\rho_{1}(d_i-\bR_i^T\balpha)-\tilde\rho_{1}(d_i)+\tilde\rho'_1(d_i)\bR_i^T\balpha-\frac 1 2 \tilde\rho_{1}''(d_i)\balpha^T\bR_i\bR_i^T\balpha\right)\right|
\leq 2 \|\balpha\|^3 \sum_{i=1}^n\|\bR_i\|^3.\end{equation*}
Note that $\sum_{i=1}^n\|\bR_i\|^3\leq h^{1/2}\sum h_i=h^{1/2}p\to0$.
Therefore, when $\|\balpha\|\leq c$, $$\E\Delta(\balpha)=\frac 1 2 \balpha^T \tilde K_n''\balpha+o(1), \text{ with }
\tilde K_n''=\sum_{i=i_0}^n \bR_i\bR_i^T \E\left[\tilde\rho''_{i,1}(d_i)\right].$$
Actually, $\tilde K_n''$ converges towards $K_0''$. Let us decompose $\|\tilde K_n''-K_0''\|$ in $$\|\tilde K_n''-K_0''\|\leq\|\tilde K_n''-K_n''\|+\|K_n''-K_0''\|.$$
The convergence to 0 of the second term is ensured by hypothesis (A3''). The first term is: $$\tilde K_n''-K_n''=\sum \bR_i^T\E(\tilde\rho''_1(d_i)-\rho''_1(d_i)).$$ The functions $\tilde \rho''_1$ and $\rho''_1$ only differ on intervals whose total length is $2/(n^2)$. Consequently, $\E(\tilde\rho''_\lambda(d_i)-\rho''_\lambda(d_i))\leq 2/(n^2) \|\tilde\rho''_\lambda-\rho''_\lambda\|_\infty \|f_\epsilon\|_\infty $ where $f_\epsilon$ denotes the density function of $\epsilon_i$. We obtain the inequality: $\|\tilde K_n''-K_n''\|\leq \frac 1 n h^{1/2} C,$ with $C$ a constant. As $h$ tends to 0 under (A4''), we deduce that $\tilde K_n''$ converges towards $K_0''$ and thus $\E\Delta(\balpha)=\frac 1 2 \balpha^T K_0''\balpha+o_\P(1).$

When $\|\balpha\|\leq c $, the convergence in probability of $|\Delta(\balpha)-\E\Delta(\balpha)|$ to 0 implies:
\begin{equation*} \left|\sum_{i=1}^n \left(\tilde\rho_{1}(d_i-\bR_i^T\balpha)-\tilde\rho_{1}(d_i)+\tilde\rho'_1(d_i)\bR_i^T\balpha\right)
-\frac 1 2 \balpha^T K_0''\balpha\right|\tend^\P 0. \end{equation*}
If $\tilde D$ and $D$ respectively denote $\tilde D:=\sum_{i=1}^n \tilde\rho_{1}(d_i-\bR_i^T\balpha)-\tilde\rho_{1}(d_i)$ and $D:=\sum_{i=1}^n \rho_{1}(d_i-\bR_i^T\balpha)-\rho_{1}(d_i)$, then $|D-\tilde D|\leq n\|\tilde\rho_{1}-\rho_{1}\|_\infty.$ Using (\ref{rho0}), we obtain the almost sure convergence of $D-\tilde D$ to 0. In the same way, if $\tilde B:=\sum_{i=1}^n\tilde\rho'_1(d_i)\bR_i^T\balpha$ and $B:=\sum_{i=1}^n\rho'_1(d_i)\bR_i^T\balpha$, we then have $|B-\tilde B|\leq n\|\tilde\rho'_{1}-\rho'_{1}\|_\infty\|\balpha\|h^{1/2}.$ When $\|\balpha\|\leq c$, properties (\ref{rho1}) implie that $B-\tilde B$ tends almost surely to 0.
All together, we have:
\begin{equation*} \left|\sum_{i=1}^n \left(\rho_{1}(d_i-\bR_i^T\balpha)-\rho_{1}(d_i)+\rho'_1(d_i)\bR_i^T\balpha\right)
-\frac 1 2 \balpha^T K_0''\balpha\right|\tend^\P 0.\\ \end{equation*}

\hspace{-0.5cm} $\bullet$ We may prove now that the convergence is uniform on the set $\{\|\balpha\|\leq c \}$.

The functions in $\balpha$: $$\sum_{i=1}^n \left(\rho_{1}(d_i-\bR_i^T\balpha)-\rho_{1}(d_i)+\rho'_1(d_i)\bR_i^T\balpha\right) ~\text{and}~ \frac 1 2 \balpha^T K_0''\balpha$$ are convex 
and the set $\{\|\balpha\|\leq c \}$ is convex, compact and independent from $n$. Proposition \ref{rock1} completes the proof.

\subsubsection*{B.4. Proof of Theorem \ref{resultat}}

\subsubsection*{B.4.1. Consistency}

In the model (\ref{modelebis}), we are willing to prove that $\hat\balpha_n=\bigcirc_\P (1).$
Let $c_n\to \infty$. We may prove that $\P\left(\|\hat\balpha_n\|>c_n\right)\to 0$.
We can deduce from (\ref{result1}) that there exists a series $c'_n$ such that $c'_n\to\infty$, $c'_n\leq c_n$ and $$\sup_{\left\{\|\balpha\|\leq c'_n \right\}} 
\left|\sum_{i=1}^n \left(\rho_{1}(d_i-\bR_i^T\balpha)-\rho_{1}(d_i)+\rho'_1(d_i)\bR_i^T\balpha\right)-\frac 1 2 \balpha^T K_0''\balpha\right|\tend^\P 0.$$ 
It is sufficient then to prove that $\P\left(\|\hat\balpha_n\|>c_n'\right)\to 0$.

\hspace{-0.5cm} $\bullet$ Suppose $\|\balpha\|=c_n'$.
 
We have $$\sum_{i=1}^n \left(\rho_{1}(d_i-\bR_i^T\balpha)-\rho_{1}(d_i)\right)=-\sum_{i=1}^n\rho'_1(d_i)\bR_i^T\balpha+\frac 1 2 \balpha^T K_0''\balpha+o_\P(1).$$ First, $$\|\frac 1 2 \balpha^T K_0''\balpha\|\geq \frac{1}{2}\underline{s}(K_0'')(c'_n)^2,$$ with $\underline{s}(K_0'')$ smallest eigenvalue of $K_0''$. As the matrix $K_0''$ is nonsingular, $\underline{s}(K_0'')>0$. Next, Proposition \ref{controle} implies that $$\|\sum_{i=1}^n\rho'_1(d_i)\bR_i^T\balpha\| =\bigcirc_\P(c_n').$$ As a consequence, the probability that the quantity \begin{equation*} \sum_{i=1}^n \rho_{1}(d_i-\bR_i^T\balpha)-\rho_{1}(d_i) = -\sum\rho'_1(d_i)\bR_i^T\balpha +\frac 1 2 \balpha^T K_0''\balpha +o_\P(1)
\end{equation*}
is negative tends to 0.
This result is true uniformly for $\balpha$ verifying $\|\balpha\|=c_n'$. We obtain: \begin{equation}\label{p1}\P\left(\inf_{\{\balpha, ~\|\balpha\|=c'_n\}} \sum_{i=1}^n \rho_{1}(d_i-\bR_i^T\balpha)-\rho_{1}(d_i) \leq 0 \right)\to 0.\end{equation}

\hspace{-0.5cm} $\bullet$ Let $\balpha$ be such that $\|\balpha\|\geq c'_n$.
 
 We define $t=\frac{c'_n}{\|\balpha\|}\in]0;1]$ and $\balpha'=t\balpha$. With the equality $d_i-\bR_i^T\balpha'=(1-t) d_i+t(d_i-\bR_i^T\balpha)$, together with the convexity of $\rho$, we have: $$\rho_{1}(d_i-\bR_i^T\balpha')-\rho_{1}(d_i)\leq t\left(\rho_{1}(d_i-\bR_i^T\balpha)-\rho_{1}(d_i)\right).$$ 
As $\|\balpha'\|=c_n'$, it comes that: \begin{equation}\label{p2}\P\left(\inf_{\{\balpha, ~\|\balpha\|\geq c'_n\}} \sum_{i=1}^n \rho_{1}(d_i-\bR_i^T\balpha)-\rho_{1}(d_i) \leq 0 \right)\to 0,\end{equation} 
or equivalently: $$\P\left(\inf_{\{\balpha, ~\|\balpha\|\geq c'_n\}} J_n(\balpha)\leq J_n(0) \right)\to 0.$$ 
The estimator $\hat\balpha_n$ has been defined as the argument realizing the minimum of $J_n$, and so, $\P\left(\|\hat\balpha_n\|\geq c'_n \right)$ tends towards zero, which achieves the proof.

\subsubsection*{B.4.2. Bahadur's representation}

We want to prove that in model (\ref{modelebis}), we have $$\hat\balpha_n={K_0''}^{-1}\left(\frac{1}{\sqrt{n}}\sum_{i=1}^n \rho_{1}'(d_i)\bR_i\right)+o_\P(1).$$

Let us first recall this result given in \shortciteA{Rockafellar}:
\begin{Proposition}
Let $\cal C$ be an open convex set. Let $f_n$ be a family of differentiable convex functions and $f$ be a differentiable convex function. If $f_n$ converges simply towards $f$ on $\cal C$, then $\nabla f_n$ converges simply towards $\nabla f$ on $\cal C$ and the convergence is uniform on every compact set of $\cal C$.
\end{Proposition}
Similarly to Proposition \ref{rock1}, this Proposition can be generalized to a convergence in probability (using Lemma \ref{ps}).

Applying this Proposition to the result (\ref{result1}) gives us that, for all $c>0$,
\begin{equation*}
\sup_{\|\balpha\|\leq c}\left|\sum_{i=1}^n \left(\rho'_{1}(d_i-\bR_i^T\balpha)\bR_i-\rho'_{1}(d_i)\bR_i\right)+ K_0''\balpha\right|\tend^\P 0.
\end{equation*}
We have proved precedently that $\hat\balpha_n=\bigcirc_\P(1)$. Then, \begin{equation}\label{bahadur}
\left|\sum_{i=1}^n \left(\rho'_{1}(d_i-\bR_i^T\hat\balpha_n)\bR_i-\rho'_{1}(d_i)\bR_i\right)+ K_0''\hat\balpha_n\right|\tend^\P 0.
\end{equation} 
By definition of $\hat\balpha_n$, $\sum_{i=1}^n \rho'_{1}(d_i-\bR_i^T\hat\balpha_n)\bR_i=0$. The convergence of (\ref{bahadur}) becomes: $$\hat\balpha_n={K_0''}^{-1}\left(\sum_{i=1}^n \rho'_{1}(d_i)\bR_i\right)+o_\P(1),$$ which is the announced result.

\subsubsection*{B.4.3. Asymptotic behavior of the functionnal part}

The model considered for this part of the proof is the model (\ref{modele}) contrarily to what precedes.

Parseval equality gives: $\|\hat f_n-f\|_2\sim\frac{1}{n}\|\hat\btheta_n-\btheta_0\|$. We decompose this bound into: $\frac{1}{n}\|\hat\btheta_n-\btheta_0\|\leq \frac{1}{n}\|\hat\btheta_n-\tilde\btheta_n\|+\frac{1}{n}\|\tilde\btheta_n-\btheta_0\|$ where $$\tilde \btheta_{i,n}= \begin{cases} z_i-\bA_i^T\bbeta_0 & \text{if~} i<i_0\\ \text{sign}(z_i-\bA_i^T\bbeta_0)\left(|z_i-\bA_i^T\bbeta_0|-\lambda\right)_+ & \text{if~} i\geq i_0 \end{cases}.$$
\shortciteA{Donoho92} proved that there exists a constant $C$ such that $\E\frac{1}{n}\|\tilde\btheta_n-\btheta_0\|\leq C \left(\frac{\log(n)}{n} \right)^{\frac{s}{1+2s}}$. The convergence in $L^2$ implies the convergence in probability. 

The term $\frac{1}{n}\|\hat\btheta_n-\tilde\btheta_n\|$ verifies the inequality $\frac{1}{n}\|\hat\btheta_n-\tilde\btheta_n\|\leq \frac{1}{n}\|A\|\|\hat\bbeta_n-\hat\bbeta_0\|+2\frac{\lambda}{n}.$ Assumptions (A2) and (A3) ensure that $\frac{1}{\sqrt{n}}\|A\|=\left(\frac{1}{n}\sum \|\bA_i\|^2 \right)^{1/2}$ is bounded and that $\|\hat\bbeta_n-\hat\bbeta_0\|=\bigcirc_\P(\frac{\lambda}{\sqrt{n}})$ through the first part of the Theorem. Then, $\frac{1}{n}\|\hat\btheta_n-\tilde\btheta_n\|=\bigcirc_\P(\frac{\lambda}{n})=\bigcirc_\P(\frac{\log(n)^{1/2}}{n}).$

\bibliographystyle{stylebib}
\bibliography{biblio}

\end{document}